\documentclass[reqno,11pt,centertags]{amsart}

\usepackage{amsmath,amsthm,amscd,amssymb,latexsym,upref,enumerate}
\input epsf
\usepackage{epsfig}
\usepackage[T2A,OT1]{fontenc}
\usepackage[ot2enc]{inputenc}
\usepackage[russian,english]{babel}
\usepackage{epic,eepicemu}
\usepackage{tikz}
\usepackage{mathrsfs}  
\usepackage{xcolor}
\usepackage[unicode,colorlinks=true,pagebackref=false]{hyperref}

\newcommand{\bbN}{{\mathbb{N}}}
\newcommand{\bbR}{{\mathbb{R}}}

\newcommand{\bbZ}{{\mathbb{Z}}}
\newcommand{\bbC}{{\mathbb{C}}}
\newcommand{\bbQ}{{\mathbb{Q}}}

\newcommand{\cU}{{\mathcal{U}}}

\newcommand{\cN}{{\mathcal{N}}}

\newcommand{\sE}{{\mathsf{E}}}
\newcommand{\ol}{\overline}

\def\bC{\overline{\mathbb{C}}}

\newcommand{\e}{{\epsilon}}

\newcommand{\z}{\zeta}

\renewcommand{\a}{\alpha}
\renewcommand{\b}{\beta}
\renewcommand{\Re}{\text{\rm Re\ }}
\renewcommand{\Im}{\text{\rm Im\ }}

\newcommand{\supp}{\operatorname{supp}}

\definecolor{Lgray}{gray}{0.7}

\allowdisplaybreaks \numberwithin{equation}{section}
\newtheorem{theorem}{Theorem}[section]
\newtheorem*{theorem-non}{Theorem}

\newtheorem{lemma}[theorem]{Lemma}
\newtheorem{proposition}[theorem]{Proposition}
\newtheorem{corollary}[theorem]{Corollary}

\theoremstyle{definition}

\newtheorem{remark}[theorem]{Remark}

\newenvironment{prf}[1]{\begin{proof}[\textbf{Proof of #1}]}{\end{proof}}

\title{Scaling limits of Chebyshev polynomials at a hard edge}

\author[B.\ Buchecker]{Benedikt Buchecker}
\address{Department of Mathematics, KU Leuven, 3001 Leuven (Heverlee), Belgium}
\email{benedikt.buchecker@kuleuven.be}

\author[B.\ Eichinger]{Benjamin Eichinger}

\address{Institute of Analysis and Scientific Computing, TU Wien, 1040 Wien,
Austria}

\email{benjamin.eichinger@tuwien.ac.uk}


\begin{document}
\begin{abstract}
    Let $\sE\subset\bbR$ be compact. We derive local scaling limits of the Chebyshev polynomials associated to $\sE$ at a gap edge. We describe the scaling limits by a one parameter family of functions explicitly given by elliptic integrals. We show that for the union of two generic intervals all functions in this family are in fact attained as limits.  
    As a consequence, we describe the local scaling behavior of the zeros of the Chebyshev polynomials at the gap edge.
\end{abstract}
	\maketitle
\section{Introduction}
In this paper we derive local scaling limits of Chebyshev polynomials at an edge. Let $\sE$ be a compact subset of $\bbC$ containing infinitely many points, and let $\|\cdot\|_\sE$ denote the sup-norm on $\sE$. The Chebyshev polynomial $T_n$ is the unique maximizer of 
\[
\sup\left\{\Re \tau_n: P(z)=\tau_nz^n+\text{l.o.t.}, \|P\|_\sE\leq 1\right\}.
\]
Chebyshev polynomials are often considered as monic polynomials minimizing the sup-norm. These extremal polynomials are easily seen to be renormalizations of $T_n$. In this work we will focus on the asymptotics of Chebyshev polynomials associated to subsets of the real line; so in the following we assume $\sE\subset\bbR$.

Studying asymptotic properties of $T_n$ is a classical topic in approximation theory going back to the work of Akhiezer, Bernstein, Faber and Szeg\H o \cite{Akhiezer1,bernstein12,Faber1920,szego:1924}. Asymptotics are described in terms of potential theoretic notions such as Green's function, capacity or equilibrium measure, see e.g. \cite{Ran95} for some background on potential theory. Root asymptotics hold for any compact $\sE\subset\bbR$ and assert convergence of $|T_n(z)|^{1/n}$. By standard arguments one can conclude from root asymptotics convergence of the normalized zero counting measure of $T_n$ to the equilibrium measure. Hence, the global asymptotic distributions of the zeros are fully understood. Comparably little is known about their local structure. This is a well-studied problem for orthogonal polynomials and typically studied through local scaling limits, so-called universality limits, of the Christoffel-Darboux kernel. For orthogonal polynomials two important universality classes are bulk universality \cite{eichinger.lukic.simanek:2026,Freud:1969,lubinsky:2009} and hard-edge universality \cite{eichinger.lukic.woracek:2024,kuijlaars.vanlessen:2002}. The first case corresponds to scaling limits given by the sine kernel and the latter one to scaling limits given by the Bessel kernel. 

In a recent work, Lubinsky \cite{lubsinksyConstr} observed that Chebyshev polynomials likewise admit local scaling limits and proved a version of bulk universality limits suitable for Chebyshev polynomials. He showed that if $I$ is an interval contained in the interior of $\sE$, then for $\xi_n\in I$ with $|T_n(\xi_n)|=1$,  one has 
\begin{align}\label{eq:lubinsky}
    \lim_{n\to\infty}\frac{T_n\left(\xi_n+\tfrac{z}{\omega_\sE(\xi_n)n}\right)}{T_n(\xi_n)}=\cos(\pi z),
\end{align}
where $\omega_\sE$ denotes the Radon-Nikodym derivative of the equilibrium measure of $\sE$ with respect to the Lebesgue measure. The goal of this paper is to study scaling limits of $T_n$ at the edge of $\sE$. In contrast to hard-edge universality limits of orthogonal polynomials and also to Lubinsky's result \eqref{eq:lubinsky}, we show that (generically) rescalings of $T_n$ at an edge do not converge to a single function, but instead the set of accumulation points is described by an explicitly given one-parameter family of functions.

Before formulating our results, we need to introduce some notation and recall some known properties of Chebyshev polynomials. In the following, let $\overline{\bbC}$ denote the Riemann sphere and $\overline\bbR=\bbR\cup\{\infty\}.$ For $\sE\subset\bbR$ compact, we call the connected components of $\overline{\bbR}\setminus\sE$, gaps of $\sE$. We denote the unbounded gap by $(\max\sE,\min\sE):=\{\infty\}\cup(\bbR\setminus [\min\sE,\max\sE])$. 

Our results rely on the geometry of $n$-extensions and on representations of extremal functions in terms of comb mappings. 
 The \textit{$n$-extension} of $\sE$ is defined by 
    \begin{align*}
            \sE_n&:=\{z\in\bbC:\ T_n(z)\in[-1,1]\}.
    \end{align*}
    The geometry of $\sE_n$ is well understood due to the alternation theorem, see e.g. \cite[Theorem 1.1]{ChriSiZi1}.    One has that $\sE_n\subset[\min\sE,\max\sE]\subset\bbR$. Moreover, if $(a,b)$ is a bounded gap of $\sE$, then $[a,b]\cap \sE_n$ can be described by two parameters $a\leq u_n\leq v_n\leq b$ so that $\sE_n\cap[a,b]=\{a,b\}\cup [u_n,v_n]$. The following extensions of the gap $(a,b)$ are possible:
    \begin{enumerate}
        \item[(i)] No extension: in this case we set $u_n=v_n=b$ (we fix this convention, as it allows us to conveniently present the results for the left gap edge.)
        \item[(ii)] Extension at one gap edge: an extension at $a$ is realized with $a=u_n<v_n\leq b$ and a similar definition at $b$
    \item[(iii)] Extension by an internal interval: this is realized with $a<u_n<v_n<b$
    \end{enumerate}
We describe scaling limits of Chebyshev polynomials at an edge by \textit{comb mappings}; cf. \cite{akhiezer.levin:1960,eremenko.yuditskii:2012,sodin.yuditskii:1992}. Our comb mappings are parametrized by the set
\begin{align*}
    \Upsilon :=\bbR \cup \{m+iy: m\in \bbZ, y\in [0,\infty]\}.
\end{align*}
Points $m+i\infty$ and $k+i\infty$ for $m\neq k$ are considered as different points. To define a topology, we glue together the points $m-0$ and $m+i\infty$ and the points $m+i0$ and $m+0$; cf. Section \ref{sec:ellipt}.
For $\tau \in \Upsilon\cap \bbC$ we define
\begin{align}\label{def:Pi}
\Pi(\tau)=\{z\in \bbC: \Im z>0, \Re z<\lfloor\Re \tau\rfloor+1\}\setminus \{\lfloor\Re\tau\rfloor+iy: 0< y\leq \Im \tau\}. 
\end{align}
In addition, we set $\Pi(m+i\infty)=\{z\in \bbC: \Im z>0, \Re z<m\}$ for $m\in\bbZ$.
The comb mapping associated to this domain is the unique conformal map $\Theta(\cdot, \tau): \bbC_+\to \Pi(\tau)$ normalized by the condition $\Theta(0,\tau)=\Re\tau$ and $\Theta(z,\tau)\sim 2i\sqrt{z}$ as $z\to +\infty$. Initially only defined in $\bbC_+$, it follows from Morera's theorem, that $\cos(\pi\Theta(z,\tau))$ extends to an entire function; cf. Proposition \ref{prop:7}. 

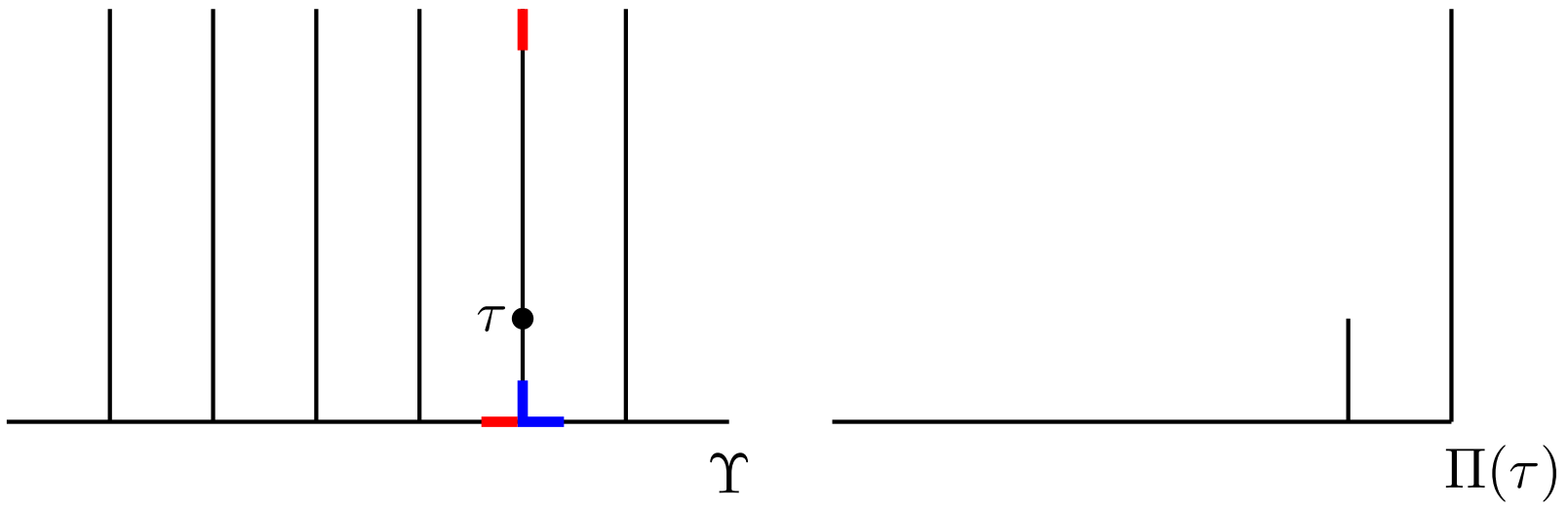
\begin{figure}
    \centering
    	\begin{tikzpicture}[scale=0.7]
		\begin{scope}[xshift=-8cm]
			\draw[-,thick] (1,0) -- (-6,0);
			\draw[-,thick] (0,0) -- (0,4);
			\draw[-,thick] (-1,0) -- (-1,4);
			\draw[-,thick] (-2,0) -- (-2,4);
			\draw[-,thick] (-3,0) -- (-3,4);
			\draw[-,thick] (-4,0) -- (-4,4);
			\draw[-,thick] (-5,0) -- (-5,4);
            \draw[-,line width=2pt,blue] (-1,0.4) -- (-1,0);
            \draw[-,line width=2pt,blue] (-1.05,0) -- (-0.6,0);
            \draw[-,line width=2pt,red] (-1.05,0) -- (-1.4,0);
            \draw[-,line width=2pt,red] (-1,4) -- (-1,3.6);
			\fill (-1,1) circle (3pt);
            \node at (-1.3,1) {$\tau$};
             \node at (1,-0.5) {$\Upsilon$};
		\end{scope}
		
		\draw[-, thick] (0,0) -- (-6,0);
		\draw[-, thick] (-1,0) -- (-1,1);
		\draw[-, thick] (0,0) -- (0,4);
			
            \node at (0.5,-0.5) {$\Pi(\tau)$};
		
	\end{tikzpicture}
    \caption{Example of a comb $\Pi(\tau)$ with $\operatorname{Im}\tau>0$. The colors indicate neighborhoods of $m$ and $m+i\infty$, respectively.}
    \label{fig:placeholder}
\end{figure}
Scaling limits of $T_n$ at a gap edge are described in terms of $\Theta(\cdot,\tau_n)$, where $\tau_n$ depends on the geometry of $\sE_n$. To define $\tau_n$ we make a case distinction. If $(a,b)$ is a bounded gap, let $u_n,v_n$ denote parameters of the $n$-extension for $(a,b)$ as defined above.  Let $\lambda_n$ be the critical point of $T_n$ in $[a,u_n]$ if $a<u_n$ and $\lambda_n=a$ if $a=u_n$. Furthermore, set $$\tau_n=\tfrac{1}{\pi}\cos^{-1}(T_n(\lambda_n))\in \Upsilon.$$ 
    If $(a,b)$ is the unbounded gap, set $\tau_n=\tfrac{1}{\pi}\cos^{-1}(T_n(a))+i\infty\in \Upsilon$.
    
    We are now ready to formulate our main theorem. We present the results for the left gap edge. With the natural adaptations in the definition of $\lambda_n,\tau_n, \Pi(\tau_n)$ and $\Theta(\cdot,\tau_n)$ all results translate also to the right gap edge. 
\begin{theorem}\label{mainTheorem}
    Let $(a,b)$ be a gap of $\sE$ and assume that there exists $\delta>0$, such that $[a-\delta,a]\subset\sE$. Let $\omega_\sE$ denote the Radon-Nikodym derivative of the equilibrium measure of $\sE$. Then the limit 
    \begin{align*}
     q:= \lim_{x\to a^-}\omega_\sE(x)\sqrt{a-x}   
    \end{align*}
exists, is positive and finite. We have
   \begin{align}\label{eq:maintheorem}
       \lim_{n\to\infty}\left(T_n\left(a+\frac{z}{(nq)^2}\right)-\cos(\pi\Theta(z,\tau_n))\right)=0
   \end{align}
   locally uniformly on $\bbC$.
\end{theorem}
\begin{remark}
\begin{enumerate}[{\rm (i)}]
    \item Since $\cos(\pi\Theta(z,\tau))=\cos(\pi\Theta(z,\tilde\tau))$ for $\tau-\tilde \tau\in 2\bbZ$, we can consider $\tau_n$ in the compact set $\Upsilon/2\bbZ$. Theorem \ref{mainTheorem} together with continuity of the map $\tau\mapsto \cos(\pi\Theta(\cdot,\tau))$ proved in Proposition \ref{prop:7} states that $\tau_{n_k}\to\hat\tau$ implies
        \begin{align}\label{eq:subseqLimit}
            \lim_{k\to\infty}T_{n_k}\left(a+\frac{z}{(n_kq)^2}\right)=\cos(\pi\Theta(z,\hat\tau)).
        \end{align}
    \item This raises the question, whether all points in $\Upsilon$ are indeed attained as subsequential limits of $\tau_{n}$. In general the answer is no. For instance, if $a=\max\sE$ all limit points must be of the form $k+i\infty$. Likewise, if $\sE=[-1,-a]\cup [a,1]$ is a union of two symmetric intervals, then all accumulation points of rescaling of $T_n$ at $-a$ will be given by $k+i\infty.$ However, we will show below that for a generic union of two intervals, the set of accumulation points is equal to the set $\Upsilon$.
    \item Limits of the form \eqref{eq:maintheorem} are commonly used to describe the asymptotics of extremal functions that do not converge but instead show asymptotically periodic or almost periodic behavior; cf.\ \cite[Theorem 12.3]{widom1969}.
    \item  Scaling limits of orthogonal polynomials at a hard-edge are typically given by a single limit. Recently, scaling limits to a whole family of limit kernels were observed for orthogonal polynomials associated to measures supported on fractal sets \cite{eichinger.lukic.kheifets.yuditskii:2026}. This type of scaling behavior reflects the self-similarity structure of the spectrum. That Chebyshev polynomials admit scaling limits as described in \eqref{eq:maintheorem}, which occur already for very simple geometries such as union of two intervals, is a new feature of Chebyshev polynomials at an edge, which is caused by the special geometry of the $n$-extensions.
\end{enumerate}
\end{remark}

For $\tau\in\bbR$, the scaling limits $\cos(\pi\Theta(z,\tau))$ at an edge are related to the scaling limits in the bulk \eqref{eq:lubinsky} by means of the transformation $z\mapsto z^2$. A similar phenomenon is known for orthogonal polynomials, see e.g. \cite[Section 10]{eichinger.lukic.woracek:2024}. 
However, for $\tau=m+ih$, with $h>0$ the scaling limits are not of this type and present new universality classes that are inherently connected to Chebyshev polynomials at an edge:

The comb mapping $\Theta(\cdot,\tau)$ for $\tau\in \bbR$ maps $\bbC_+$ onto a quarter plane and is given by 
\[
    \Theta(z,\tau) = \lfloor \tau \rfloor +1 +2i\sqrt{z-v},\quad v=\left(\frac{1+\lfloor \tau \rfloor -\tau}{2}\right)^2.
\]
Similarly, for $\tau=m+i\infty$, we have
\[\Theta(z,m+i\infty) = m+ 2i\sqrt{z},\]
which can be viewed as the natural analogue of \eqref{eq:lubinsky} for gap edges. However, for $\tau=m+ih$ with $m\in \bbZ,h\in(0,\infty),$ $\Theta(\cdot, \tau)$ corresponds to a new universality class given by elliptic integrals. More precisely, for such $\tau$ there exist  $0<\lambda<u<v<\infty$ such that the comb mapping is given by 
\[
\Theta(z,m+ih)=m+ i\int_0^z \frac{\zeta-\lambda}{\sqrt{\zeta(\zeta-u)(\zeta-v)}}d\zeta,
\]
where the branch of the square root is chosen so that $\sqrt{\zeta(\zeta-u)(\zeta-v)}>0$ for $\zeta>v$. Analyzing the behavior of this Christoffel-Schwarz mapping shows that it maps $\bbC_+$ to a comb domain $\Pi(\tau)$ and that $\Theta^{-1}((-\infty,m+1],\tau)=(-\infty,0]\cup[u,v]$.
The degenerate cases can be well understood geometrically in terms of the interval $[u,v]$. The case $ \tau=m+i0$ for $m\in \bbZ$ corresponds to the degeneration $u\to 0$ and the case $\tau =k+i\infty$ corresponds to the case $u\to\infty$.

To formulate corollaries about the local structure of zeros of $T_n$ we will label the zeros relative to $a$ by $x_k^{(n)}$ so that
\[
\dots<x_{-2}^{(n)}<x_{-1}^{(n)}<x_{0}^{(n)}\leq a.
\]
When we write a limit statement involving some $x_{k}^{(n)}$, it is implied as part of that statement that $x_{k}^{(n)}$ exists for all large enough n. We will use the notation $f \lesssim g $  if $\limsup f /g \leq C$ for some $C<\infty$ and $f \asymp  g$ if $f\lesssim g$ and $g \lesssim f$.

\begin{corollary}\label{cor:1}
    Let $(a,b)$ be a gap of $\sE$ and assume that there exists $\delta>0$, such that $[a-\delta,a]\subset\sE$. Then for any $j<l$, with $j,l\in\bbZ_{\leq 0}$ 
    \begin{align}\label{eq:zeros}
    x_{l}^{(n)}-x_j^{(n)}\asymp \frac{1}{n^2}.
    \end{align}
    Moreover, for any $R>0$ and large enough $n$, $T_n$ can have at most one zero in $[a,a+Rn^{-2}]$, and along a subsequence where such a zero $x_1^{(n_j)}$ exists, \eqref{eq:zeros} holds also for $j<l=1$.
\end{corollary}

In Theorem \ref{mainTheorem} we have exhausted all cases of limits which might appear as edge scaling limits of Chebyshev polynomials. However, it is not clear whether all these limits can actually appear. In particular, given $\tau=m+ih$, for $h\in(0,\infty),$ does there exist $\sE$ and a subsequence $n_k$ so that rescalings of $T_{n_k}$ converge to $\cos(\pi\Theta(\cdot,\tau))$? We give an affirmative answer to this question for the simplest non-trivial case, i.e., when $\sE$ is a union of two intervals, which was already studied by Akhiezer \cite{Akhiezer1,Akhiezer2,Akhiezer3}. Strikingly, a much stronger statement is true: For generic $(a;b)$, the Chebyshev polynomial for $\sE=[-1,a]\cup [b,1]$ rescaled at $a$ admits all possible limit points $\cos(\pi\Theta(\cdot,\tau))$ for $\tau\in \Upsilon$. To indicate the dependence on $\sE$, we write $T_n^\sE$ for the Chebyshev polynomials associated to the compact $\sE$.

\begin{theorem}\label{thm:all_limits}
    There exists a dense $G_\delta$ set $F\subseteq\{(a;b): -1 <a<b<1\}$, i.e. $F$ is a countable intersection of dense open sets, such that for all $\sE=[-1,a]\cup [b,1]$ with $(a;b)\in F$ and $\tau \in \Upsilon$ there exists a subsequence $n_k$ with $T_{n_k}^\sE\left(a+\frac{z}{(n_kq)^2}\right)\xrightarrow{k\to\infty}\cos (\pi\Theta(z,\tau))$ locally uniformly for $z\in\bbC$. 
\end{theorem}

Our proofs rely on the fact that Chebyshev polynomials $T_n$ admit a representation 
\begin{align*}
   T_n(z)=\cos(\pi n\theta_n(z)), 
\end{align*}
in terms of a comb mapping $\theta_n$ associated to an $n$-regular comb $\Pi_n$. Statements such as \eqref{eq:subseqLimit} can intuitively be understood by the geometric picture that rescalings of the comb $\Pi_n$ converge to $\Pi(\hat{\tau})$. These geometric ideas can be formalized using Carath\'eodory's kernel convergence \cite{Pom92}. We do not follow this method in this paper, but instead use more elementary techniques based on Herglotz functions and Christoffel-Schwarz mappings. 

In Section 2 we recall the necessary prerequisites on Chebyshev polynomials, Herglotz functions and potential theory. In Section 3 we derive the representations of the comb mappings $\Theta(\cdot,\tau)$ as elliptic integrals and prove continuous dependence on the parameter $\tau$. In Section 4 we prove Theorem \ref{mainTheorem} and Corollary \ref{cor:1} and in Section 5 we prove Theorem \ref{thm:all_limits}.

\section{Preliminary}
This section is mainly of preliminary nature. We recall basic facts about Herglotz functions, comb representations of Chebyshev polynomials and prove some preparatory lemmas.

\subsection{Logarithmic potential theory and comb mappings}
Logarithmic potential theory plays a pivotal role in describing the asymptotics of Chebyshev polynomials. For subsets of the real line, explicit formulas for Chebyshev polynomials in terms of Green's functions are available. All of what we present below is a consequence of the well-known alternation theorem, see e.g. \cite[Theorem 1.1]{ChriSiZi1}. 

Let $\sE\subset\bbR$ compact and of positive capacity $\operatorname{Cap}(\sE)>0$. We call $\xi_0<\xi_1<\dots<\xi_{M_P}$, $\xi_j\in\sE$, a set of alternation points for a polynomial $P$ with $\|P\|_\sE\leq 1$, if $P(\xi_j)=(-1)^{M_P-j}$. The alternation theorem asserts that the Chebyshev polynomial $T_n$ is uniquely characterized by the property that it has a maximal set of alternation points, i.e., $M_{T_n}=n$. As a direct consequence, we get that the $n$-extensions $\sE_n = T_n^{-1}([-1,1])$ is contained in $[\min\sE,\max\sE]$.  We call $\Pi_n$ an $n$-regular comb, if it is of the form 
	\[
	\Pi_n=\left\{x+iy: \frac{-k_-}{n}<x<\frac{k_+}{n}, y>0\right\}\setminus \bigcup_{k=1}^{n-1}\left\{\frac{-k_-+k}{n}+iy:\ 0\leq y\leq h_k\right\},
	\]
 	for $k_-\in\bbN,k_+\in 2\bbN$, $k_-+k_+=n$ and $h_k \geq 0$.

    The boundary of $\Pi_n$ is locally connected, and therefore any conformal map from $\bbC_+$ to $\Pi_n$ can be continuously extended to the boundary. Therefore, there exists a unique conformal map $\theta_n:\bbC_+\to\Pi_n$ normalized so that $\theta_n(\min\sE)=\tfrac{-k_-}{n},  \theta_n(\max\sE)=\tfrac{k_+}{n}$ and $\theta_n(\infty)=\infty$. 

    For any $\sE\subset\bbR$ compact, there exists an $n$-regular comb domain $\Pi_n$ and an associated mapping $\theta_n$ so that $\sE_n=\theta_n^{-1}([\tfrac{-k_-}{n},\tfrac{-k_+}{n}])$ and 
	\begin{align}\label{sec2:combRepCheby}
	T_n(z)=\cos(\pi n\theta_n(z)),
	\end{align}
	see e.g. \cite[Section 5.4]{sodin.yuditskii:1992} or \cite[proof of Theorem 3.1]{Ahlfors}. Note that although $\theta_n$ is initially only defined on $\bbC_+,$ it can be shown that the right-hand side of \eqref{sec2:combRepCheby} extends to a meromorphic function on $\overline{\bbC}$. Since $T_n(\max \sE)=1$, the map $\theta_n$ (and comb $\Pi_n$) is unique up to adding a constant of the form $2\frac{k}{n}$ for $k\in \bbZ$. Therefore, for any $a\in\sE$, we might assume $n\theta_n(a)\in [0,2)$.

    Moreover, the equilibrium measure of $\sE_n$ is absolutely continuous and given by $\mu_{\sE_n}((u,v))=\Re \theta_n(v)-\Re\theta_n(u)$ and the Green's function is given by $g_{\sE_n}(z)=\pi\Im \theta_n(z)$. 

    It follows from root asymptotics of $T_n$ and \eqref{sec2:combRepCheby}, that the Green's function and the equilibrium measure of $\sE_n$ converge to the ones of $\sE.$

     \begin{lemma}\label{lem:14}
        For $\sE_n = T_n^{-1}([-1,1])$, it holds that
        \[
        g_{\sE_n}(z) \xrightarrow{n\to \infty}g_{\sE}(z)
        \]
        locally uniformly for $z\in \bbC\setminus [\min \sE,\max \sE]$.
    \end{lemma}
    \begin{proof}
        We use \eqref{sec2:combRepCheby}, i.e., 
        \[
        T_n(z) = \cos(\pi n\theta_n(z))=\frac{1}{2}\Big(e^{in\pi\theta_n(z)}+ e^{-in\pi\theta_n(z)}\Big).
        \]
        Using  $g_{\sE_n}(z)=\pi\Im \theta_n(z)$, this implies 
        \begin{equation}\label{eq:root_ineq}
        1 - e^{-2ng_{\sE_n}(z)}\leq \frac{2|T_n(z)|}{e^{ng_{\sE_n}(z)}}\leq 1+e^{-2ng_{\sE_n}(z)}.
        \end{equation}
        Theorem 3.2 in \cite{ChriSiZi1} establishes root asymptotics, i.e., 
        \[
        |T_n(z)|^{1/n}\xrightarrow{n\to \infty} e^{g_{\sE}(z)}
        \]
        locally uniformly for $z\in \bbC\setminus [\min\sE,\max\sE]$. Taking the $n$-th root in \eqref{eq:root_ineq}, combined with the observation that for any compact $K\subset\bbC\setminus [\min \sE,\max\sE]$ there exists $C>0$ such that for $z\in K$
        \[
        g_{\sE_n}(z) \geq g_{[\min\sE,\max\sE]}(z)\geq C>0,
        \]
        proves the assertion.
    \end{proof}

It is standard, see e.g. \cite{Simon:2007}, that convergence of the Green's functions implies convergence of the equilibrium measures.

\begin{lemma}\label{lem:13}
        Let $\sE\subseteq K_n \subseteq [\min\sE,\max\sE]$ be a sequence of compact sets with $g_{K_n}(z)\to g_{\sE}(z)$ for all $z\in U$, where $U\subseteq \bbC\setminus [\min\sE,\max\sE]$ is nonempty and open. Then,
        \[
        \mu_{K_n}\xrightarrow{n\to \infty}\mu_{\sE}
        \]
        in the weak$^*$ topology. In particular, $\mu_{\sE_n}\to \mu_\sE$ for the $n$-extensions $\sE_n$. 
    \end{lemma}
    \begin{proof}
        Let us first note that the convergence $g_{K_n}\to g_\sE$ holds locally uniformly in $\bbC\setminus[\min\sE,\max\sE]$. In this domain, both $g_{K_n}$ and $g_\sE$ are harmonic and $g_{\sE}-g_{K_n}\geq 0$ because $K_n\supseteq \sE$. It is a consequence of Harnack's theorem (see e.g. \cite[Theorem 1.3.10]{Ran95}) that this means that either $g_\sE-g_{K_n}\to \infty$ locally uniformly or there exists a convergent subsequence $g_{\sE}-g_{K_{n_k}} \to h$ where $h$ is harmonic in $\bbC\setminus [\min\sE,\max\sE]$. Since $g_{\sE}-g_{K_n}\to 0$ in $U$, we know that the first case is impossible and furthermore $h\equiv 0$ by the identity principle. This implies $g_{K_n}\to g_\sE$ locally uniformly in $\bbC\setminus[\min\sE,\max\sE]$.
        
        Now, let us prove the convergence of the measures $\mu_{K_n}$. Since $\mu_{K_n}$ are probability measures, for every subsequence we can find a further subsequence which is convergent, i.e. $\mu_{K_{n_k}}\to \mu$ with $\mu$ a probability measure supported on $[\min \sE,\max\sE]$. We consider the logarithmic energy 
        \begin{align*}
         I(\nu) := \int_\bbC\int_\bbC \log |z-w|d\nu(z)d\nu(w) 
     \end{align*}
     of a compactly supported measure $\nu$. By definition of the equilibrium measure
     \[
     I(\mu_{\sE})\leq I(\mu_{K_{n}}) \leq I(\mu_{[\min\sE,\max\sE]}),
     \]
     which means that by taking a further subsequence we can suppose $I(\mu_{K_{n_k}})\to C\in \bbR$. Then, for $z\in\bbC\setminus [\min\sE,\max\sE]$
     \[
     g_{K_{n_k}}(z)=\int_{K_{n_k}}\log |z-w|d\mu_{K_{n_k}}(w) -I(\mu_{n_k})\xrightarrow{k\to \infty} \int_\bbR \log |z-w|d\mu(w) -C.
     \]
     On the other hand,
     \[
     g_{K_{n_k}}(z) \xrightarrow{k\to\infty}g_{\sE}(z)= \int_{\sE}\log |z-w|d\mu_{\sE}(w) -I(\mu_{\sE})
     \]
     which implies
     \[
      \int_{\sE}\log |z-w|d\mu_{\sE}(w)= \int_\bbR \log |z-w|d\mu(w) +I(\mu_{\sE})-C \leq  \int_\bbR \log |z-w|d\mu(w).
     \]
     By \cite[Theorem A.21]{Simon:2007}, we know $\mu=\mu_{\sE}$ which proves the lemma.
    \end{proof}

\subsection{Herglotz functions}
Let $\mu$ be a probability measure on $\bbR$ and $\Phi_\mu(z)=-\int\log(z-t)d\mu(t)$, then 
\[
\Phi_\mu'(z)=\int\frac{d\mu(t)}{t-z},\quad z\in\bbC_+.
\]
Hence, $\Phi_\mu'$ is a Herglotz function, that is an analytic map from $\bbC_+$ to $\bbC_+$. This observation is used for instance in the study of ergodic Schr\"odinger operators \cite{kotani:2015} and was used in the deterministic setting in \cite{eichinger.lukic:2025}. In our setting, we have
	\[
	\pi i\theta_n'(z)=\int\frac{d\mu_{\sE_n}(t)}{t-z}=:m_{\sE_n}(z),
	\]
	where $\mu_{\sE_n}$ is the equilibrium measure for $\sE_n$, see \cite[Section 4]{eremenko.yuditskii:2012}. We will first study the scaling behavior of $m_{\sE_n}$ and then derive asymptotic statements for $\theta_n$. Combining this with \eqref{sec2:combRepCheby}, will yield Theorem \ref{mainTheorem}.

We call $f$ a Herglotz function if it maps $\bbC_+$ analytically into $\overline{\bbC_+}=\bbC_+\cup\bbR\cup\{\infty\}$ and denote this class by $\overline{\cN}$.  
If $f\not\equiv \infty$, then there exists $\alpha\geq 0,\beta\in\bbR$ and a Poisson finite measure $\mu$ so that 
	\begin{align}\label{sec2:eq1}
	    f(z)=\alpha z+\beta+\int_{-\infty}^\infty\left(\frac{1}{t-z}-\frac{t}{1+t^2}\right)d\mu(t),\quad \int\frac{d\mu(t)}{1+t^2}<\infty.
	\end{align}
    Using the chordal metric on $\overline{\bbC_+}$, we  equip $\overline{\cN}$ with the metric of local uniform convergence in $\bbC_+$, which makes it a  compact metric space. If $f_n,f\in\overline{\cN}\setminus \{\infty\}$ and $\a_n,\a,\b_n,\b,\mu_n,\mu$ are the parameters in the respective integral representations \eqref{sec2:eq1}, then the Grommer--Hamburger theorem \cite[Proposition 7.8]{lukic:book} asserts that $f_n\to f$ if and only if 
    \begin{align*}
		\a_n+\int&\frac{d\mu_n(t)}{1+t^2}\to \a+\int\frac{d\mu(t)}{1+t^2},\quad  \b_n\to \b,\\
	    &\mu_n\to\mu,\quad \text{in the $w^*$-topology of $C_c(\bbR)$.}
	\end{align*}

Our asymptotic analysis of $m_{\sE_n}$ essentially relies on the following two lemmas:
\begin{lemma}\label{lem:5}
Let $f_n, f\in\overline{\cN}\setminus \{\infty\}$ and assume that $f_n\to f$. Let $\mu_n,\mu$ be the measures in the integral representation of $f_n$ and $f$, respectively. If there exists an interval $I$ and $n_0>0$ such that  for $n\geq n_0$ $\supp(\mu_n)\cap I = \emptyset$. Then, $f_n$ and $f$ have holomorphic extensions to $\bbC_+\cup\bbC_-\cup I$ and $f_n$ converges locally uniformly to $f$ on this domain.
\end{lemma}

\begin{proof}
It's convenient to rewrite \eqref{sec2:eq1} to 
\begin{align}\label{sec2:eq2}
   f_n(z)=\a_nz+\b_n+\int\frac{1+tz}{t-z}\frac{d\mu_n(t)}{1+t^2} 
\end{align}
and similarly for $f$. Since 
\[
f_n(i)=\b_n+i\left(\a_n+\int_\bbR\frac{d\mu_n(t)}{1+t^2}\right),
\]
we conclude that the sequences $\b_n,\a_n+\int_\bbR\tfrac{d\mu_n(t)}{1+t^2}$ are bounded. Furthermore, a simple calculation proves the following inequality. If $1\leq R<\infty$, $z\in \bbC$ with $\text{dist}(z,\supp \mu_n)>R^{-1}$ and $|z|\leq R$, then
    \begin{align*}
        \left|\frac{1+zt}{t-z}\right|\leq 4R^3, \quad t\in \supp \mu_n.
    \end{align*}
    Now, it is clear that \eqref{sec2:eq2} defines a holomorphic function on $\bbC_+\cup\bbC_-\cup I$ with $f_n(z)= \overline{ f_n(\ol z)}$. The same holds for $f$ since the weak$^*$ convergence $\mu_n\to \mu$ implies $\supp(\mu)\cap I=\emptyset$. Furthermore, the sequence $f_n$ is uniformly bounded on compact subsets of $\bbC_+\cup\bbC_-\cup I$ and forms a normal family. Since $f_n\to f$ in $\bbC_+$  we conclude that $f_n\to f$ locally uniformly in $\bbC_+\cup\bbC_-\cup I$. 
\end{proof}
In the following, we will extend Herglotz functions $f:\bbC_+\to\bbC_+$ to $\bbC_+\cup \bbR$ by setting $f(x)=\lim_{\varepsilon\to 0^+} f(x+ i\varepsilon)$ for $x\in \bbR$ if this limit exists.
\begin{lemma}\label{lem:6}
    Let $\sE\subseteq K_n\subset [\min\sE,\max\sE]$ be a sequence of compact sets and $(a,b)$ a bounded gap of $\sE$ and assume that there is $\delta>0$ so that $[a-\delta,a]\subset\sE$. Assume that
    \begin{align*}
        K_n \cap [a,b]= \{a,b\} \cup [u_n,v_n],
    \end{align*}
    with $a\leq u_n\leq v_n\leq b$. Let $\lambda_n$ and $\Tilde{\lambda}_n$ be the points where $g_{K_n}$ attains its maximum in the intervals $[a,u_n]$ and $[v_n,b]$ respectively. If $\lim_{n\to\infty}\mu_{K_n}=\mu_\sE$, then $\lim_{n\to\infty}(v_n-u_n)=0$ and 
    \begin{align*}
        m_{K_n}(x) \frac{\sqrt{(x-a)(x-u_n)(x-v_n)}}{(x-\lambda_n)(x-\Tilde{\lambda}_n)} \xrightarrow{n\to\infty} m_{\sE}(x) \frac{\sqrt{x-a}}{x-\lambda}
    \end{align*}
    uniformly for $x$ in compact subsets of $\bbC_+\cup(a-\delta,b)$, where the maximum of $g_\sE$ in $[a,b]$ is attained in $\lambda$.
\end{lemma}
\begin{remark}
\begin{enumerate}
    \item[(i)]
    Note that the interval $[a,u_n]$ and/or $[v_n,b]$ may collapse to a single point where we define $\lambda_n=a=u_n$ and/or $\Tilde{\lambda}_n=v_n=b$. Otherwise the open interval $(a,u_n)$ or $(v_n,b)$ is a part of the complement $\bbR\setminus K_n$. In $\bbR\setminus K_n$ we know $g_{K_n}''=-m_{K_n}' <0$ which means that $g_{K_n}$ is strictly concave and admits a unique maximum in $[a,u_n]$ or $[v_n,b]$. Similarly, $g_\sE$ admits a unique maximum in $[a,b]$.
    \item[(ii)] For the formulation of the lemma, we have assumed that the gap $(a,b)$ is bounded. However, the same proof shows that
    \[
    m_{K_n}(x)\sqrt{x-a}\xrightarrow{n\to \infty}m_\sE(x)\sqrt{x-a},
    \]
    uniformly in compact subsets of $\bbC_+\cup (a-\delta,\infty)$ if $(a,b)$ is unbounded.
\end{enumerate}
\end{remark}
\begin{proof}
Let us first prove that $v_n-u_n\to 0.$ Assume that there exists a subsequence, so that $\lim_{j\to\infty}(v_{n_j}-u_{n_j})=c>0$. By passing to a further subsequence, we can assume that $v_{n_j}\to v$, and $u_{n_j}\to u$ so that $v-u=c.$ Then, for large enough $j$, $[u_{n_j},v_{n_j}]\supseteq [u+c/4,v-c/4]$ and $\mu_{K_{n_j}}([u+c/4,v-c/4])\geq \mu_{[\min\sE,\max\sE]}([u+c/4,v-c/4])>0$, by \cite[Theorem A.15]{Simon:2007}. The Portmanteau theorem implies $$0<\limsup_{j}\mu_{K_{n_j}}([u+c/4,v-c/4])\leq\mu_\sE([u+c/4,v-c/4])$$. On the other hand,  since $[u+c/4,v-c/4]\cap\sE=\emptyset$ we have $\mu_\sE([u+c/4,v-c/4])=0$. A contradiction.

    In the interval $(a-\delta,b)$ the Weyl-function $m_{K_n}$ has the following behavior. On $(a-\delta,a)$ it is purely imaginary with imaginary part $\pi\omega_{K_n}$ with singularity at $a$. If $u_n>a$, then in $(a,u_n)$ the function $m_{K_n}$ is real, starting from $-\infty$ at $a$, changing sign at $\lambda_n\in (a,u_n)$ and going to $+\infty$ at $u_n$. In $(u_n,v_n)$, $m_{K_n}$ is again purely imaginary with imaginary part $\pi\omega_{K_n}$. If $v_n<b$, then in $(v_n,b)$ it is real and changes sign at $\Tilde{\lambda}_n$.

    Thus, we may define the function
    \begin{align*}
        \eta_n: (a-\delta,b)\to \bbR, t\mapsto \begin{cases}
            \frac{1}{2}, &t \leq a,\\
            1, &a< t\leq \lambda_n,\\
            0, &  \lambda_n< t\leq u_n,\\
            \frac{1}{2}, &u_n<t \leq v_n,\\
            1, &v_n< t\leq \Tilde{\lambda}_n,\\
            0, &  \Tilde{\lambda}_n< t.
        \end{cases}
    \end{align*}
    Then, $\arg m_{K_n} = \pi \eta_n(t)$ and
    \begin{align*}
        \log m_{K_n}(z) - \int_{\bbR} \frac{1}{t-z}\eta_n(t)dt
    \end{align*}
    defines a Herglotz function which extends to $(a-\delta,b)$ as a real function. Furthermore, if we apply the exponential function, this is still a Herglotz function which is positive on $(a-\delta,b)$. We obtain
    \begin{align*}
        f_n(z) := m_{K_n}(z)\frac{\sqrt{(z-a)(z-u_n)(z-v_n)}}{(z-\lambda_n)(z-\Tilde{\lambda}_n)}.
    \end{align*}
    Let us take any limit point of this sequence of Herglotz functions, i.e. $f_{n_k}\to f$. By Lemma \ref{lem:5} this convergence holds uniformly on compact subsets of $\bbC_+\cup\bbC_-\cup(a-\delta,b)$. By taking a further subsequence we can suppose that $u_{n_k}, v_{n_k}\to u \in [a,b]$, $\lambda_{n_k}\to \rho\in [a,b]$ and $\Tilde{\lambda}_{n_k} \to \Tilde{\rho}\in [a,b]$. Since $\mu_{K_n}\to\mu_\sE$, we know $m_{K_n}\to m_\sE$ in $\bbC_+$. Then, for $z\in \bbC_+\cup (a-\delta,b)$ we have
    \begin{align*}
        f(z) = m_\sE(z) \frac{\sqrt{(z-a)(z-u)^2}}{(z-\rho)(z-\Tilde{\rho})}
    \end{align*}
    Since $f_{n_k}$ is positive on $(a-\delta,b)$, the same holds for $f$. In particular, $f$ has no zeros or poles. The function $m_\sE$ has a square root singularity at $a$ and potentially a zero at $\lambda$ (in particular if $\lambda<b$). Then we know $\{\rho,\Tilde{\rho}\}=\{u,\lambda\}$ and
    \begin{align*}
        f(z) = m_\sE(z) \frac{\sqrt{z-a}}{(z-\lambda)}.
    \end{align*}
    This proves the assertion.
\end{proof}

\section{Properties of elliptic integrals}\label{sec:ellipt}
In this section we study properties of the comb mappings $\Theta(\cdot,\tau)$ defined  by \eqref{def:Pi} and the subsequent paragraph.  We derive explicit formulas as Christoffel-Schwarz mappings and based on those, we show continuous dependence on $\tau.$

We start with deriving explicit formulas for $\Theta(\cdot,\tau)$. For $\tau=m+iy$ with $m\in\bbZ$ and $y>0$, it can be represented as a Christoffel-Schwarz mapping. Let us point out that the shape of the comb $\Pi(\tau)$ enforces natural restrictions on the parameters of this mapping. Firstly, the base of the comb $\Pi(\tau)$ is at the same height to the left and to the right of the needle. Secondly, the needle is placed at distance one to the ``imaginary boundary'' of the comb. In terms of $\Theta$, these conditions translate to  $\Theta(0,\tau)=\Theta(u,\tau)$ and $\Theta(v,\tau)-\Theta(u,\tau)=1$.  It will be convenient to express them in terms of elliptic integrals. Recall that the complete elliptic integral of first and second kind with modulus $k$ are given by
\begin{align*}
    K(k):= \int_0^{\pi/2} \frac{1}{\sqrt{1-k^2\sin^2\theta}}d\theta,\\
    E(k):= \int_0^{\pi/2} \sqrt{1-k^2\sin^2\theta}d\theta.
\end{align*}
\begin{lemma}\label{lem:17}
    If $\tau\in\bbR,$ we have
   \[
    \Theta(z,\tau) = \lfloor \tau \rfloor +1 +2i\sqrt{z-v},\quad \text{where }v=\left(\frac{1+\lfloor \tau \rfloor -\tau}{2}\right)^2.
\]
If $\tau=m+i\infty$ for $m\in\bbZ$, then 
\[\Theta(z,m+i\infty) = m+ 2i\sqrt{z}.\] 
If $\tau=m+iy$ with $m\in\bbZ$ and $y>0$, then there exist $0<\lambda<u<v$ so that 
\begin{align}\label{eq:CSformulaTheta}
    \Theta(z,\tau) = m + i \int_0^z \frac{\zeta-\lambda}{\sqrt{\zeta(\zeta-u)(\zeta-v)}}d\zeta.
\end{align}
The parameters $\lambda,u,v$ are subject to the conditions
\begin{align}\label{eq:ell_int_1}
0&=(v-\lambda)K\left(\sqrt{\frac{u}{v}}\right)- vE\left(\sqrt{\frac{u}{v}}\right),\\
\frac{1}{2}&=\sqrt{u}E\left(i\sqrt{\frac{v}{u}-1}\right)-\frac{\lambda}{\sqrt{u}} K\left(i\sqrt{\frac{v}{u}-1}\right).\label{eq:ell_int_2}
\end{align}
Moreover, for any $m\in\bbZ$ and $0<\lambda<u<v$ satisfying \eqref{eq:ell_int_1} and \eqref{eq:ell_int_2} there exists $y>0$ so that for $\tau=m+iy$ the mapping $\Theta(\cdot,\tau)$ is of the form \eqref{eq:CSformulaTheta}.
\end{lemma}
\begin{proof}
   We start with $\tau \in \bbR$. In this case $\Pi(\tau)$ is a quarter plane and the only vertex is $\lfloor \tau \rfloor +1$ with angle $\frac{\pi}{2}$. Let $\Theta(v,\tau)=\lfloor \tau\rfloor +1$, then $\Theta(z,\tau)\sim 2i\sqrt{z}$ as $z\to \infty$ implies
\[
    \Theta(z,\tau) = \lfloor \tau \rfloor +1 +2i\sqrt{z-v}.
\]
From $\Theta(0,\tau)=\tau$, we get $v= \left(\frac{1+\lfloor \tau \rfloor -\tau}{2}\right)^2$.

If $\tau=m+i\infty$, we also have a quarter plane. Same as before, we get for $m\in \bbZ$
\[\Theta(z,m+i\infty) = m+ 2i\sqrt{z}.\] 

If $\tau =m+iy$ with $m\in \bbZ$ and $y\in (0,\infty)$, then we have four vertices $m-$, $\tau$, $m+$ and $m+1$ ($m$ has to be counted twice since we have a slit). We already know $\Theta(0,\tau)=m-$ and if we define $0<\lambda<u<v$ by 
\[
\Theta(\lambda,\tau)=\tau,\quad \Theta(u,\tau) = m+,\quad \Theta(v,\tau)=m+1,
\]
we get
\[
\Theta(z,\tau) = m + i \int_0^z \frac{\zeta-\lambda}{\sqrt{\zeta(\zeta-u)(\zeta-v)}}d\zeta.
\]
We know 
\[
0=\Theta(u,\tau) -\Theta(0,\tau) =\int_0^u \frac{\zeta-\lambda}{\sqrt{\zeta(\zeta-u)(\zeta-v)}}d\zeta.
\]
Substituting $\zeta = u\sin^2\theta$ and simplifying shows \eqref{eq:ell_int_1}.
Similarly, $\Theta(v,\tau)-\Theta(u,\tau)=1$ implies \eqref{eq:ell_int_2}.

Conversely, given $m\in\bbZ$ and $0<\lambda<u<v$, we can form the Christoffel-Schwarz mapping \eqref{eq:CSformulaTheta}, which maps the upper half plane into a polygonal region with fixed angles at the vertices $0,\lambda, u, v$. If this mapping is denoted by $\Psi$, the conditions \eqref{eq:ell_int_1} and \eqref{eq:ell_int_2} imply that $\Psi(0)=\Psi(u)$ and $\Psi(v)-\Psi(u)=1$. Thus, with $y=\operatorname{Im}\Psi(\lambda)$, we get that $\Psi=\Theta(\cdot,m+iy).$
\end{proof}

Next, we want to show that $\Theta(\cdot,\tau)$ depends continuously on $\tau\in \Upsilon$. The first step is to consider convergent sequences of $\lambda_k,u_k,v_k$ that satisfy \eqref{eq:ell_int_1} and \eqref{eq:ell_int_2}. 
\begin{lemma}\label{lem:18}
    Let $0<\lambda_k<u_k<v_k<\infty$ satisfy \eqref{eq:ell_int_1} and \eqref{eq:ell_int_2} for every $k$. Furthermore, suppose that $\lim_{k\to \infty}\lambda_k=\lambda_0, \lim_{k\to \infty}u_k=u_0$ and $\lim_{k\to \infty}v_k=v_0$ where $\lambda_0,u_0,v_0\in [0,\infty]$. Then, one of the following holds:
    \begin{enumerate}
        \item[(i)]$0=\lambda_0=u_0<v_0<\infty$.
        \item[(ii)]  $0<\lambda_0<u_0<v_0<\infty$ and they satisfy \eqref{eq:ell_int_1} and \eqref{eq:ell_int_2}. 
        \item[(iii)]  $\lambda_0=u_0=v_0=\infty$ and $\lim_{k\to\infty} \frac{\lambda_k}{u_k}=\lim_{k\to\infty} \frac{u_k}{v_k}=1$. 
    \end{enumerate}
\end{lemma}
\begin{proof}
    We start by proving that in any case we have $v_0>0$. Indeed, we know that for all $k$
\[
1=\int_{u_k}^{v_k}\frac{\zeta-\lambda_k}{\sqrt{\zeta(\zeta-u_k)(v_k-\zeta)}}d\zeta.
\]
On the other hand, for $0<\lambda_k<u_k\leq\zeta\leq v_k$, we have
\[
\frac{\zeta-\lambda_k}{\sqrt{\zeta}}\leq \sqrt\z\leq\sqrt{v_k}
\]
and 
\[
\int_{u_k}^{v_k}\frac{1}{\sqrt{(\zeta-u_k)(v_k-\zeta)}}d\zeta=\pi.
\]
We conclude that 
\[
\int_{u_k}^{v_k}\frac{\zeta-\lambda_k}{\sqrt{\zeta(\zeta-u_k)(v_k-\zeta)}}d\zeta\leq \pi\sqrt{v_k}.
\]
Hence, we get $1/\pi^2\leq v_k$ and consequently $1/\pi^2\leq v_0$. 

We continue by showing that $v_0=\infty$ implies $\lambda_0=u_0=\infty$. Let us assume that $\lambda_0<\infty$. From \eqref{eq:ell_int_1} it follows that
    \begin{align*}
        \frac{E\left(\sqrt{\frac{u_k}{v_k}}\right)}{K\left(\sqrt{\frac{u_k}{v_k}}\right)}\xrightarrow{k\to\infty} 1-\frac{\lambda_0}{v_0}=1.
    \end{align*}
    This implies $\frac{u_k}{v_k}\to 0$ and
    \begin{align*}
        K\left(i\sqrt{\frac{v_k}{u_k}-1}\right)\xrightarrow{k\to\infty}0.
    \end{align*}
    Using that $\lambda_k<u_k$, this also implies $\frac{\lambda_k}{\sqrt{u_k}}K\left(i\sqrt{\frac{v_k}{u_k}-1}\right)\to 0$. Then, from \eqref{eq:ell_int_2} using that $E(ik)/k\geq 1$ for $k>0$ we get
    \begin{align*}
        \max (\sqrt{u_k},\sqrt{v_k-u_k}) \leq \sqrt{u_k}E\left(i\sqrt{\frac{v_k}{u_k}-1}\right)\to \frac{1}{2}.
    \end{align*}
    This is a contradiction, since either $u_k\to \infty$ or $v_k-u_k\to \infty$. Thus, we have shown that if one of $\lambda_0,u_0,v_0$ is $\infty$, then all of them are $\infty$. Moreover, \eqref{eq:ell_int_2} shows that if $\lambda_0=u_0=v_0=\infty$, then
    \begin{align*}
        0&\leq E\left(i\sqrt{\frac{v_k}{u_k}-1}\right) - K\left(i\sqrt{\frac{v_k}{u_k}-1}\right)\\
        &\leq E\left(i\sqrt{\frac{v_k}{u_k}-1}\right) -\frac{\lambda_k}{u_k}K\left(i\sqrt{\frac{v_k}{u_k}-1}\right) \xrightarrow{k\to \infty} 0.
    \end{align*}
    This implies $\lim_{k\to \infty}\frac{u_k}{v_k}=\lim_{k\to\infty}\frac{\lambda_k}{u_k}=1$ in this case. 

   It remains to show that if $0<u_0<\infty$, then $0<\lambda_0<u_0<v_0<\infty$. We already know $0\leq \lambda_0\leq u_0\leq v_0<\infty$. If $u_0=v_0$, then $K\left(\sqrt{\frac{u_k}{v_k}}\right)\xrightarrow{k\to \infty} \infty$. Then \eqref{eq:ell_int_1} implies that $\lambda_0=v_0=u_0$ and
    \[
    \sqrt{u_k}E\left(i\sqrt{\frac{v_k}{u_k}-1}\right)-\frac{\lambda_k}{\sqrt{u_k}} K\left(i\sqrt{\frac{v_k}{u_k}-1}\right) \xrightarrow{k\to \infty} \sqrt{u_0}\frac{\pi}{2} - \sqrt{u_0}\frac{\pi}{2}=0.
    \]
    However, by \eqref{eq:ell_int_2} the left hand side is $1/2$ which is a contradiction. Thus, $u_0<v_0$. By continuity, we get that $\lambda_0,u_0,v_0$ satisfy \eqref{eq:ell_int_1} and \eqref{eq:ell_int_2}. The first is equivalent to 
    \[
    \int_0^{u_0} \frac{\zeta-\lambda_0}{\sqrt{\zeta(\zeta-u_0)(\zeta-v_0)}}d\zeta=0.
    \]
    If either $\lambda_0=0$ or $\lambda_0=u_0$, then there is no sign change in the integral. However, the integrand is nonzero which would be a contradiction. Therefore, if $0<u_0<\infty$, then $0<\lambda_0<u_0<v_0<\infty$ and they satisfy \eqref{eq:ell_int_1} and \eqref{eq:ell_int_2}. 
\end{proof}

Our goal is to show that $\Theta(\cdot,\tau)$ depends continuously on $\tau$. In order for this to be true, we need to define the correct topology on $\Upsilon$. Observe that for $\tau_\varepsilon = m-\varepsilon$ one has $ \lfloor \tau_\varepsilon\rfloor= m-1$ and uniformly for $z$ in compact subsets of $\bbC_+\cup \bbR\setminus \{0\}$
\[
\begin{aligned}
\Theta(z,\tau_\varepsilon)&= \lfloor \tau_\varepsilon\rfloor +1 + 2i\sqrt{z-\left(\frac{\lfloor \tau_\varepsilon\rfloor +1-\tau_\varepsilon}{2}\right)^2}= m+2i \sqrt{z-\frac{\varepsilon^2}{4}}\\
&\xrightarrow{\varepsilon\to 0}m+2i\sqrt{z}=\Theta(z,m+i\infty)\ne \Theta(z,m).
\end{aligned}
\]
Accordingly, we will equip $\Upsilon =\bbR \cup \{m+iy: m\in \bbZ, y\in (0,\infty]\}$ with the usual topology of $\bbC$ except that $\lim_{\varepsilon\to 0^+}(m-\varepsilon)=m+i\infty= \lim_{y\to \infty}(m+iy)$. 

To account for this we need to change the neighborhood bases of $m$ and $m+i\infty$; cf Figure \ref{fig:placeholder}. To define the topology at $m$ we define for $a,b\in (0,1)$ the sets $\cU_{a,b}=\{m+iy\mid y\in [0,a)\}\cup [m,m+b)$. Then $\{\cU_{a,b}\mid a,b\in(0,1)\}$ form a neighbourhood base at $m$. Likewise, we define for $\e,\delta\in (0,1)$ the set $\mathcal V_{\e,\delta}:=(m-\e,m)\cup\{m+iy\mid y>1/\delta\}\cup \{m+i\infty\}$. A neighbourhood base at $m+i\infty$ is given by $\{\mathcal V_{\e,\delta}\mid\ \e,\delta\in(0,1)\}.$ With this modification, we get the following statement.

\begin{proposition}\label{lem:11}
    Let $\tau_k\in \Upsilon$ with $\tau_k\to \tau \in \Upsilon$. Then, $\Theta(z,\tau_k)\to\Theta(z,\tau)$ locally uniformly for $z\in \bbC_+$.    
\end{proposition}
\begin{proof}
    Let us start with $\tau_k\in \bbR$. If $\tau_k\to \tau$ with $\tau\in \bbR$, then $\lfloor\tau_k\rfloor=\lfloor\tau\rfloor$ for $k$ large enough. We get
    \[
    \begin{aligned}
        \Theta(z,\tau_k)&= \lfloor \tau_k\rfloor +1 + 2i\sqrt{z-\left(\frac{\lfloor \tau_k\rfloor +1-\tau_k}{2}\right)^2}\\
        &\xrightarrow{k\to \infty}\lfloor \tau\rfloor+1+2i\sqrt{z-\left(\frac{\lfloor \tau\rfloor +1-\tau}{2}\right)^2}=\Theta(z,\tau),
    \end{aligned}
\]
locally uniformly for $z\in \bbC_+$. Let us note that $\Theta(\cdot,\tau_k)$ and $\Theta(\cdot,\tau)$ extend continuously to $\bbR$ and the convergence also holds uniformly on compact subsets of $ \bbR\setminus\{v\}$ where $v=\left(\frac{\lfloor \tau\rfloor +1-\tau}{2}\right)^2$. If $\tau_k \to m+i\infty$ with $m\in\bbZ$, the assertion follows from the discussion before the proposition. 

Now, let $\tau_k = m+iy_k$ with $(0,\infty)\ni y_k\to y\in [0,\infty]$ and
    \begin{align*}
        \Theta(z,\tau_k) &= m + i \int_0^z \frac{\zeta-\lambda_k}{\sqrt{\zeta(\zeta-u_k)(\zeta-v_k)}}d\zeta,\\
        \Theta(z,\tau) &= m + i \int_0^z \frac{\zeta-\lambda}{\sqrt{\zeta(\zeta-u)(\zeta-v)}}d\zeta.
    \end{align*}
Since $0<\lambda_k<u_k<v_k$, we know that for $z\in \bbC_+$ with $\Re z\leq\lambda_k$ and $\zeta\in [0,z]$
\begin{equation}\label{eq:sqrt_bound}
\Big|\frac{\zeta-\lambda_k}{\sqrt{\zeta(\zeta-u_k)(\zeta-v_k)}}\Big| \leq \frac{1}{\sqrt{|\zeta|}},
\end{equation}
which implies $|\Theta(z,\tau_k)|\leq |m|+ 2\sqrt{|z|}$. From the fact that $\Theta(\cdot,\tau_k)$ are Herglotz functions, this already implies that subsequentially $\Theta(\cdot,\tau_k)$ converges locally uniformly in $\bbC_+$ (to a Herglotz function that is not constant $\infty$). Thus, by taking a subsequence (denoted the same) we can assume $\Theta(\cdot,\tau_k)\to \Theta_0$. By taking a further subsequence, we may also suppose $\lambda_k \to \lambda_0,u_k\to u_0,v_k\to v_0$, where $0\leq \lambda_0\leq u_0\leq v_0\leq \infty$. 
    
    By Lemma \ref{lem:18}, there are three cases for the limits $\lambda_0,u_0,v_0$. If we denote them as case $(i)$, $(ii)$ and $(iii)$, then
    \[
    \Theta'(\zeta,\tau_k) = i \frac{\zeta-\lambda_k}{\sqrt{\zeta(\zeta-u_k)(\zeta-v_k)}} \xrightarrow{k\to \infty} \Theta_0 '(\zeta)= \begin{cases}
            i\frac{1}{\sqrt{\zeta-v_0}},& \text{case }(i),\\
            i\frac{\zeta-\lambda_0}{\sqrt{\zeta(\zeta-u_0)(\zeta-v_0)}} , & \text{case }(ii), \\
            i\frac{1}{\sqrt{\zeta}}, &\text{case }(iii),
        \end{cases}
    \]
    uniformly for $\zeta$ in compact subsets of $\bbC_+\cup \bbR\setminus \{0,u_0,v_0\}$. In particular, $\Theta_0$ extends to $\bbR$ and the convergence $\Theta(\cdot,\tau_k)\to \Theta_0$ is uniform on compact subsets of $\bbC_+\cup \bbR\setminus \{0,u_0,v_0\}$. Furthermore, for $z\in \bbC_+$ with $\Re z<\lambda_0$, using \eqref{eq:sqrt_bound} we get 
    \[
        \Theta(z,\tau_k)\xrightarrow{k\to \infty}m+ i \int_0 ^z \Theta_0'(\zeta)d\zeta.
    \]
    This implies for $z\in \bbC_+\cup \bbR$
    \[
    \Theta_0(z)=\begin{cases}
            m+2\sqrt{v_0} +2i\sqrt{z-v_0},& \text{case }(i),\\
            m+i\int_0^z\frac{\zeta-\lambda_0}{\sqrt{\zeta(\zeta-u_0)(\zeta-v_0)}}d\zeta , & \text{case }(ii), \\
            m+2i\sqrt{z}, &\text{case }(iii).
        \end{cases}
    \]
    
    In case $(ii)$, we know that $\lambda_0,u_0,v_0$ satisfy \eqref{eq:ell_int_1} and \eqref{eq:ell_int_2}. By Lemma \ref{lem:17} we have $\Theta_0=\Theta(\cdot,m+i\Tilde{y})$ for some $\Tilde{y}>0$. Since $0<\lambda_0<u_0$, it holds that \[m+i\Tilde{y}=\Theta(\lambda_0,m+i\Tilde{y})=\Theta_0(\lambda_0)=\lim_{k\to\infty}\Theta(\lambda_k,\tau_k)=\lim_{k\to\infty}\tau_k=\tau.\] 
    Thus, $\Theta(\cdot,\tau_k)\to \Theta(\cdot,\tau)$.
    
    In case $(iii)$, for every $M>0$ we have $\lim_{k\to\infty}\Theta(M^2,\tau_k)=\Theta_0(M^2)=m+2i M$. Since $M^2<\lambda_k$ for $k$ large enough, we have $M<\Im \tau_k$ for $k$ large enough. Thus, $\tau=\lim_{k\to \infty}\tau_k =m+i\infty$ and
    \[
    \lim_{k\to\infty}\Theta(\cdot,\tau_k)= \Theta_0 = \Theta(\cdot,\tau).
    \]
    
    In case $(i)$ consider any $x>v_0$. Then, $x>v_k$ for $k$ large enough. Then, $\Re\Theta(x,\tau_k)=m+1$. On the other hand
    \[
    \Theta(x,\tau_k)\xrightarrow{k\to\infty} \Theta_0(x) = m+2\sqrt{v_0} +2i\sqrt{x-v_0}. 
    \]
    Thus, $v_0 = 1/4$ and $\Theta_0= \Theta(\cdot,m)$. Furthermore, from \eqref{eq:sqrt_bound} we also get
    \[
    \Im \tau_k =\Im \Theta(\lambda_k,\tau_k)\leq 2\sqrt{\lambda_k},
    \]
    which implies $\tau=\lim_{k\to \infty}\tau_k =m$. 
    
    Thus, $\Theta(\cdot,\tau_k)\to \Theta(\cdot,\tau)$ in every case.
\end{proof}

\begin{proposition}\label{prop:7}
    Let $\tau\in \Upsilon$. Then, the function $\cos(\pi \Theta(\cdot,\tau))$ extends to an entire function. Furthermore, if $\Upsilon\ni\tau_k\to \tau\in \Upsilon$, then $\cos(\pi \Theta(\cdot,\tau_k))$ converges to $\cos(\pi \Theta(\cdot,\tau))$ locally uniformly in $\bbC$. 
\end{proposition}
\begin{proof}
    We can extend $\cos(\pi\Theta(\cdot,\tau))$ to $\bbC_-$ by setting 
    \[
    \cos(\pi\Theta(z,\tau)) = \overline{\cos(\pi\Theta(\overline{z},\tau))}, \quad z\in \bbC_-.
    \]
    We know that $\Theta(\cdot,\tau)$ extends to $\bbR$ and that for $x\in\bbR$ either $\Theta(x,\tau)\in \bbR$ or $\Theta(x,\tau)=m+iy$, where $m\in \bbZ,y>0$. In both cases $\cos(\pi \Theta(x,\tau))\in \bbR$ which makes the extension continuous on $\bbC$. By Morera's theorem it is an entire function. In the proof of Proposition \ref{lem:11} we have seen that the convergence $\Theta(\cdot,\tau_{k})\to \Theta(\cdot,\tau)$ is uniform on compact subsets of $\bbC_+\cup \bbR$ with at most four points removed. Then, Cauchy's theorem implies that $\cos(\pi\Theta(\cdot,\tau_k))\to \cos(\pi\Theta(\cdot,\tau))$ locally uniformly in $\bbC$. 
\end{proof}
In Corollary \ref{cor:1}, we will prove statements about the zeros of Chebyshev polynomials. We will require some properties of the zeros of the limit functions $\cos(\pi \Theta(\cdot,\tau))$ collected in the following lemma.
\begin{lemma}\label{lem:15}
    For every $\tau\in \Upsilon$, $\cos(\pi \Theta(\cdot,\tau))$ has infinitely many simple zeros on the negative real line and at most one simple zero on the positive real line. Let them be ordered by
    \[
    \dots< x_{-3}^{\tau} < x_{-2}^{\tau} <x_{-1}^{\tau}< x_{0}^{\tau} \leq 0 \,(<x_1^\tau).
    \]
    Then, for every $j<l\leq 0$
    \begin{equation}\label{eq:sep_zeros_Theta}
    0<\inf_{\tau \in \Upsilon} (x_l^\tau -x_j^\tau) \leq \sup_{\tau \in \Upsilon} (x_l^\tau -x_j^\tau)<\infty.
    \end{equation}
    The inequalities hold for $j<l=1$ if one takes the infimum and supremum over those $\tau$ for which $x_1^\tau$ exists.
\end{lemma}
\begin{proof}
    We know that
    \begin{equation}\label{eq:zeros_theta}
    \{z\in\bbC: \cos(\pi\Theta(z,\tau))=0\} = \Theta(\cdot,\tau)^{-1}\Big(\Big\{m+\frac{1}{2}: m\in \bbZ\Big\}\Big).
    \end{equation}
    In a small neighborhood $N$ around each such point we have that $\Theta(\cdot,\tau)$ satisfies $\Im\Theta(z,\tau)>0$ for $z\in N\cap \bbC_+$ and $\Theta(x,\tau)\in \bbR$ for $x\in N\cap \bbR$. Thus, $\Theta'(x,\tau)\neq0$ for every $x\in N\cap\bbR$. Since the cosine has simple zeros, this implies that all zeros of $\cos(\pi \Theta(\cdot,\tau))$ are simple too.
    
    Since $ \Theta(\bbR_-,\tau)=(-\infty,\Re\tau)$, \eqref{eq:zeros_theta} implies that $\cos (\pi \Theta(\cdot,\tau))$ has infinitely many zeros on the negative real axis $\bbR_-$. The only other possible zero is the preimage of $\Re \tau +\frac{1}{2}$. However, $\Re\tau+\frac{1}{2}\in \overline{\Pi(\tau)}$ if and only if $\tau\in \bbR$ with $\tau-\lfloor\tau\rfloor <\frac{1}{2}$ or $\Im \tau \in (0,\infty)$. Therefore, these are the cases where $x_1^\tau$ exists. Since $\Theta(\cdot,\tau_{k})\to \Theta(\cdot,\tau)$ if $\tau_k \to \tau$, Hurwitz's theorem also shows $x_l^{\tau_k}\to x_l^\tau$ (or $x_l^{\tau_k}\to x_{l-1}^\tau$ if $\tau_k \to (m+\frac{1}{2})^-$ for $m\in \bbZ$). Since $\cos(\pi \Theta(\cdot,\tau+2m))=\cos (\pi\Theta(\cdot,\tau))$ for $m\in \bbZ$ and $\{\tau\in \Upsilon:0\leq\Re \tau\leq 2\}$ is compact, this implies \eqref{eq:sep_zeros_Theta}.
\end{proof}

\section{Scaling limits}
In this section we prove normality and study subsequential limits of 
\begin{align}\label{eq:3Nov}
T_n\left(a+\frac{z}{(nq)^2}\right),\quad z\in\bbC,
\end{align}
leading to the proof of  Theorem \ref{mainTheorem}. Throughout the section, let $\sE\subset\bbR$ be compact and let $(a,b)$ be a gap of $\sE$ and assume that there exists $\delta>0$, such that $[a-\delta,a]\subset\sE$.

 Let $\omega_\sE$ denote the Radon-Nikodym derivative of the equilibrium measure $\mu_\sE$ with respect to the Lebesgue measure. By \cite[Theorem A.13]{Simon:2007}, the limit
 $$q= \lim_{x\to a^-}\omega_\sE(x)\sqrt{a-x}>0$$
 exists and is positive and finite. 

 By \eqref{sec2:combRepCheby}, we have $T_n(z) =\cos (\pi n\theta_n(z))$, where $\Im \theta_n =\frac{1}{\pi}g_{\sE_n}$ and $\sE_n=T_n^{-1}([-1,1])$. Thus, our task reduces to proving normality of the family \eqref{eq:3Nov} and to study subsequential limits of 
\[
n\theta_n\left(a+ \frac{z}{(nq)^2}\right),\quad z\in\bbC_+.
\]
  
The proposition below is formulated for general $K_n\supseteq \sE$ and conformal maps $\theta_{K_n}$ given by $\Im \theta_{K_n} =\frac{1}{\pi}g_{K_n}$. Note that we still have $\pi i\theta_{K_n}'=m_{K_n}$. We will apply Proposition \ref{prop:1} for $K_n\equiv E$ and $K_n=\sE_n.$ In particular, for the choice $K_n\equiv\sE$ we get normality of the family \eqref{eq:3Nov}. 
\begin{proposition}\label{prop:1}
    Let $\sE\subseteq K_n\subset \bbR$ be a sequence of compact sets and $(a,b)$ a gap of $\sE$ and assume that there is $\delta>0$ so that $[a-\delta,a]\subset\sE$. Assume that
    \begin{align*}
        K_n \cap [a,b]= \{a,b\} \cup [u_n,v_n],
    \end{align*}
    with $a\leq u_n\leq v_n\leq b$, as well as $\mu_{K_n}\to\mu_\sE$ as $n\to\infty$. Let $K_{n_k}$ be a subsequence of $K_n$. Then the following holds
    \begin{enumerate}
        \item[(i)] If there exists $\delta>0$ such that $|u_{n_k}-a|>\delta$ for all $k$ or $(a,b)$ is unbounded, then for every $R>0$, there exists a sequence $\varepsilon_k\to 0$ as $k\to\infty$ such that 
    \begin{align*}
        n_kg_{K_{n_k}}\left(a+\frac{z}{(n_kq)^2}\right) \leq 2\pi(1+\varepsilon_k)\sqrt{|z|}, \quad |z|\leq R.
    \end{align*}
    Furthermore,
    \begin{align*}
        n_k\theta_{K_{n_k}}\left(a+\frac{z}{(n_kq)^2}\right)-n_k\theta_{K_{n_k}}(a) \xrightarrow{k\to\infty}2 i\sqrt{z},
    \end{align*}
    uniformly for $z$ in compact subsets of $\bbC_+\cup \bbR$.
    \end{enumerate}
    \item[(ii)] If for every $k$ it holds that $u_{n_k}=a$, then
    \begin{align*}
        n_k\theta_{K_{n_k}}\left(v_{n_k}+\frac{z}{(n_kq)^2}\right)-n_k\theta_{K_{n_k}}(v_{n_k}) \xrightarrow{k\to\infty}2 i\sqrt{z},
    \end{align*}
    uniformly for $z$ in compact subsets of $\bbC_+\cup \bbR$.
 \end{proposition}
\begin{proof}
    Let us start with case (i). The assumptions, combined with Lemma \ref{lem:6} show
    \begin{align*}
        m_{K_{n_k}}(z)\sqrt{z-a}\xrightarrow{k\to\infty} m_E(z) \sqrt{z-a},
    \end{align*}
    locally uniformly in a neighborhood of $a$. For $z\in \bbC_+\cup \bbR$ let
    \begin{align*}
        G_k(z) :=n_k\theta_{K_{n_k}}\left(a+\frac{z}{(n_kq)^2}\right)-n_k\theta_{K_{n_k}}(a) .
    \end{align*}
    Let $R>0$. For big enough $k$, $G_k$ has an extension to the interval $[-R,R]$. Then, for $z\in \bbC_+\cup [-R,R]$  
    \begin{align*}
        G_k(z) &= \int_0^z G_k'(\zeta)d\zeta = \int_0^z \frac{-in_k}{\pi(n_kq)^2}m_{K_{n_k}}\left(a+\frac{\zeta}{(n_kq)^2}\right) d\zeta\\
        &=\int_0^z \frac{-in_k}{(n_kq)^2\pi}m_{K_{n_k}}\left(a+\frac{\zeta}{(n_kq)^2}\right) \frac{\sqrt{\frac{\zeta}{(n_kq)^2}}}{\sqrt{\frac{\zeta}{(n_kq)^2}}}d\zeta.
    \end{align*}
    If we define
    \begin{align*}
        \Tilde{q}_k(\zeta):=\frac{1}{\pi q}m_{K_{n_k}}\left(a+\frac{\zeta}{(n_kq)^2}\right) \sqrt{\frac{\zeta}{(n_kq)^2}},
    \end{align*}
    then
    \begin{align*}
        G_k(z) = \int_0^z \Tilde{q}_k(\zeta) \frac{-i}{\sqrt{\zeta}} d\zeta.
    \end{align*}
    Now, uniformly for $\zeta$ in compact subsets of $\bbC_+\cup \bbR$ we have 
    \begin{align*}
        \Tilde{q}_k(\zeta)\xrightarrow{k\to\infty} \frac{1}{\pi q}\lim_{z\to a}m_\sE(z)\sqrt{z-a} = -\frac{1}{q}\lim_{t\to a^-}\omega_\sE(t)\sqrt{a-t}=- 1,
    \end{align*}
    where the penultimate equality follows from \cite[Theorem A.13]{Simon:2007} together with Perron--Stieltjes inversion.
    Thus, with
    \begin{align*}
        \varepsilon_k:= \max_{|\zeta|\leq R} |\Tilde{q}_k(\zeta)+ 1|\xrightarrow{k\to\infty} 0,
    \end{align*}
    we have for $|z|\leq R$,
    \begin{align*}
        |G_k(z)| \leq \int_0^{z} |\Tilde{q}_k(\zeta)|\frac{1}{\sqrt{|\zeta|}}|d\zeta| \leq (1+\varepsilon_k) \int_0^{|z|} \frac{1}{\sqrt{t}}dt= 2(1+\varepsilon_k)\sqrt{|z|}.
    \end{align*}
    This implies 
    \begin{align*}
        n_kg_{K_{n_k}}\left(a+\frac{z}{(n_kq)^2}\right) = \pi\Im G_k(z) \leq 2\pi(1+\varepsilon_k)\sqrt{|z|}
    \end{align*}
    Since $g_{K_{n_k}}(\overline{z}) = g_{K_{n_k}}(z)$, this bound extends to $z\in\bbC_-$ as well. 
    Now, uniformly for $z$ in compact subsets of $\bbC_+\cup \bbR$, we have
    \begin{align*}
        G_k(z) = \int_0^z \Tilde{q}_k(\zeta) \frac{-i}{\sqrt{\zeta}} d\zeta\xrightarrow{k\to\infty} \int_0^z \frac{ i}{\sqrt{\zeta}}d\zeta = 2i\sqrt{z},
    \end{align*}
    which proves (i).
    The case (ii) follows in the same way from 
    \[
    m_{K_{n_k}}(z)\sqrt{z-v_{n_k}}\xrightarrow{k\to\infty} m_E(z) \sqrt{z-a}.
    \]
\end{proof}
As an immediate consequence we get normality of the family \eqref{eq:3Nov}.
\begin{corollary}\label{cor:3}
    For every $R>0$, there exist $\varepsilon_n>0$ with $\varepsilon_n \to 0$ such that 
    \[
    \left|T_n\Big(a+\frac{z}{(nq)^2}\Big)\right|\leq \exp\left(2\pi(1+\varepsilon_n)\sqrt{|z|}\right), \quad |z|\leq R.
    \]
    In particular, the rescaled polynomials $z\mapsto T_n\Big(a+\frac{z}{(nq)^2}\Big)$ form a normal family.
\end{corollary}
\begin{proof}
    Case (i) in Proposition \ref{prop:1} with $K_n =\sE$ and $u_n=v_n=b$ shows that for every $R>0$, there exist $\varepsilon_n\to 0$ with
    \begin{align*}
        ng_{\sE}\left(a+\frac{z}{(nq)^2}\right) \leq 2\pi (1+\varepsilon_n)\sqrt{|z|}, \quad |z|\leq R.
    \end{align*}
    Then, the first assertion follows from the Bernstein-Walsh inequality
    \[
    \left|T_n\Big(a+\frac{z}{(nq)^2}\Big)\right| \leq \exp\left(ng_{\sE}\left(a+\frac{z}{(nq)^2}\right) \right) \|T_n\|_\sE
    \]
    and $\|T_n\|_\sE \leq 1$. Montel's theorem proves the normality.
\end{proof}
This already implies that subsequential limits exist. Our goal is to show that they are exactly the limits given in Theorem \ref{mainTheorem}. Let $[u_n,v_n]$ denote the additional interval of $E_n$ in the gap $(a,b)$, $\lambda_n$ the extremal value of $g_{\sE_n}$ in $[a,u_n]$ and $\tau_n = \frac{1}{\pi}\cos^{-1}(T_n(\lambda_n))=n\theta_n(\lambda_n)\in \Upsilon$. Since we can shift $\theta_n$ by $2\frac{k}{n}$ for $k\in \bbZ$, we can assume that $n\theta_n(a)=\Re \tau_n\in [0,2)$. That means for every subsequence, we find a further subsequence $n_k$ such that $\tau_{n_k}\to \tau\in \Upsilon$. Thus, in light of Proposition \ref{prop:7}, we have to show that in this case
\[
T_{n_k}\left(a+\frac{z}{(n_kq)^2}\right)\xrightarrow{k\to\infty}\cos(\pi\Theta(z,\tau)).
\]
By normality, it is enough to prove this for $z\in \bbC_+$. By taking further subsequences, we can suppose that we are in one of three cases: Either $|u_{n_k}-a|>\delta$ for some $\delta>0$, $u_{n_k}=a$ for all $k$, or $u_{n_k}>a$ and $u_{n_k}\to a$. Let us analyze them separately.
\subsection{Positive gap}
\begin{proposition}\label{prop:6}
    Let $n_k$ be a subsequence such that $\tau_{n_k}\to \tau\in \Upsilon$ and $|u_{n_k}-a|>\delta$ for some $\delta>0$ or $(a,b)$ is unbounded. Then, $\tau=m+i\infty$ for $m\in\bbZ$ and 
    \[
    T_{n_k}\left(a+\frac{z}{(n_kq)^2}\right)\xrightarrow{k\to\infty}\cos\left(\pi\Theta(z,\tau)\right),
    \]
    locally uniformly for $z\in \bbC$.
\end{proposition}
\begin{proof}
If $|u_{n_k}-a|>\delta$ for some $\delta>0$, then \[\Im \tau_{n_k} =\frac{1}{\pi}n_k g_{\sE_{n_k}}(\lambda_n)\to \infty.\] Moreover, $\Re\tau_n =n_k\theta_{n_k}(a)\in \{0,1\}$. Then, $\tau_{n_k}\to \tau\in \Upsilon$ means that $\Im \tau=\infty$, and $n_k\theta_{n_k}(a)=\Re\tau_{n_k}=\Re \tau$ for $k$ large enough. By Lemma \ref{lem:13}, $\sE_n=T_n^{-1}([-1,1])$ fulfills the premise of Proposition \ref{prop:1} and we get
\begin{align*}
    T_{n_k}\left(a+\frac{z}{(n_kq)^2}\right)&=\cos\left(\pi n_k\theta_{n_k}\left(a+\frac{z}{(n_kq)^2}\right)\right)\\
    &\xrightarrow{k\to \infty} \cos(\pi\Re \tau+ 2\pi i \sqrt{z})=\cos\left(\pi\Theta(z,\tau)\right)
\end{align*}
for $z\in \bbC_+$. The convergence extends to locally uniform convergence in $\bbC$ by Corollary \ref{cor:3} and Proposition \ref{prop:7}.
\end{proof}

\subsection{No gap}
\begin{proposition}\label{prop:5}
    Let $n_k$ be a subsequence such that $\tau_{n_k}\to \tau\in \Upsilon$ and $u_{n_k}=a$ for all $k$. Then, $\tau\in\bbR$ or $\tau=m+i\infty$ and
    \[
    T_{n_k}\left(a+\frac{z}{(n_kq)^2}\right)\xrightarrow{k\to\infty}\cos\left(\pi\Theta(z,\tau)\right),
    \]
    locally uniformly for $z\in \bbC$.
\end{proposition}
\begin{proof}
Since $u_{n_k}=a$, we have $\lambda_{n_k}=a$, $v_{n_k}>u_{n_k}$ and $\tau_{n_k}\in \bbR$. We further know that \[ n_k\theta_{n_k}(v_{n_k})= \left\lfloor
n_k\theta_{n_k}(a)\right\rfloor+1 =\lfloor \tau_{n_k}\rfloor +1.
\]
Since $\bbR\ni\tau_{n_k}\to \tau\in \Upsilon$, the convergence in $\Upsilon$ shows that for $k$ large enough $n_k\theta_{n_k}(v_{n_k})= m$, where $m=\lfloor \tau\rfloor+1$ if $\tau \in \bbR$ and $m=\Re \tau$ if $\Im \tau =\infty$. Proposition \ref{prop:1} (ii) shows that
\begin{equation}\label{eq:shifted_theta_n}
 n_k\theta_{n_k}\left(v_{n_k}+\frac{z}{(n_kq)^2}\right) \xrightarrow{k\to\infty}m+2 i\sqrt{z},
\end{equation}
uniformly for $z$ in compact subsets of $\bbC_+\cup \bbR$ and
\begin{equation}\label{eq:shifted_rescaling} 
    T_{n_k}\left(v_{n_k}+\frac{z}{(n_kq)^2}\right)\xrightarrow{k\to \infty} \cos \left(m\pi+2 \pi i\sqrt{z}\right),
\end{equation}
locally uniformly for $z\in \bbC$. We want to rescale at $a$ which means changing $z$ to $z-(n_kq)^2(v_{n_k}-a)$.

Let $\xi_k$ be the next alternation point of $T_{n_k}$ smaller or equal than $a$ or equivalently $\xi_k\leq a$ such that $v_{n_k}-\xi_k$ is minimal and $T_{n_k}(\xi_k)=(-1)^{m-1}$. Then, Hurwitz's theorem combined with \eqref{eq:shifted_rescaling} shows that $(n_kq)^2(v_{n_k}-\xi_k) \to \frac{1}{4}$. Thus, $(n_kq)^2(v_{n_k}-a) \leq (n_kq)^2(v_{n_k}-\xi_k)$ is bounded. By taking a further subsequence we can suppose that $(n_kq)^2(v_{n_k}-a) \to v$. From \eqref{eq:shifted_theta_n}, we get
\[
 \tau_{n_k}=n_k\theta_{n_k}\left(a\right) = n_k\theta_{n_k}\left(v_{n_k}+\frac{(n_kq)^2(a-v_{n_k})}{(n_kq)^2}\right) \xrightarrow{k\to \infty} m -2\sqrt{v}.
\]
Then, knowing that $\tau_{n_k}\to \tau$, we get $v= \left(\frac{\lfloor \tau\rfloor +1-\tau}{2}\right)^2$ if $\tau \in \bbR$ and $v=0$ if $\Im \tau =\infty$. Furthermore, from \eqref{eq:shifted_rescaling} we get
\begin{align*}
T_{n_k}\left(a+\frac{z}{(n_kq)^2}\right)&= T_{n_k}\left(v_{n_k}+\frac{z-(n_kq)^2(v_{n_k}-a)}{(n_kq)^2}\right)\\
&\xrightarrow{k\to \infty} \cos \left(m\pi+2\pi i\sqrt{z-v}\right)= \cos\left(\pi\Theta(z,\tau)\right).
\end{align*}
\end{proof}

\subsection{Shrinking gap}
We will consider subsequences ${n_k}$ where $\tau_{n_k}\to \tau\in \Upsilon$,  $u_{n_k}>a$ and $u_{n_k}\to a$ as $k\to\infty$. The analysis of this section is much more delicate than the previous ones and requires a further case distinction. By assumption, we know that $\Im \tau_{n_k}=n_k\Im \theta_{n_k}(\lambda_{n_k})\to y_0\in[0,\infty].$ We will consider the cases $y_0=0, y_0=\infty$ and $y_0\in(0,\infty)$ separately. In particular, the latter case gives rise to the elliptic integrals discussed in the Section 3.

As before, we define 
\begin{equation}\label{eq:def_g_k}
        G_k(z) :=n_k\theta_{n_k}\left(a+\frac{z}{(n_kq)^2}\right),
\end{equation}
for $z\in \bbC_+\cup \bbR$. Furthermore, we introduce the auxiliary variables
\begin{align}
    \Tilde{\lambda}_k :=(n_kq)^2(\lambda_{n_k}-a),\\
    \Tilde{u}_k :=(n_kq)^2(u_{n_k}-a),\\
    \Tilde{v}_k :=(n_kq)^2(v_{n_k}-a).\label{eq:tildeV}
\end{align}
We will investigate the behavior of the comb map $G_k$. We know that since $E_n\supseteq E$, Proposition \ref{prop:1} shows
\begin{align*}
    \limsup_{k\to\infty}\text{Im}(G_k(z)) &=\frac{1}{\pi} \limsup_{k\to\infty}n_kg_{\sE_{n_k}}\Big(a+ \frac{z}{(n_kq)^2}\Big)\\
    &\leq \frac{1}{\pi}\limsup_{k\to\infty}n_kg_{\sE}\Big(a+ \frac{z}{(n_kq)^2}\Big) \leq 2\sqrt{|z|}.
\end{align*}
In particular,
\begin{align}\label{eq:bound_imag_nu}
      \limsup_{k\to\infty}\Im\tau_{n_k} = \limsup_{k\to\infty}(\Im G_k(\Tilde{\lambda}_k))\leq 2 \limsup_{k\to\infty} \sqrt{\Tilde{\lambda}_k}.
\end{align}
First, we will prove a convenient formula for $G_k$. This will also build the bridge to the elliptic integrals from Section 3. 
\begin{lemma}
    Let $G_k$ be defined by \eqref{eq:def_g_k}. For $z\in \bbC_+\cup \bbR$, it holds that
    \begin{align}\label{eq_christoffel_schwarz}
        G_k(z) =\Re\tau_{n_k}+ i\int_0^z \Tilde{h}_k(\zeta) \frac{\zeta- \Tilde{\lambda}_k}{\sqrt{\zeta(\zeta-\Tilde{u}_k)(\zeta-\Tilde{v}_k)}}d\zeta,
    \end{align}
    where  $\Tilde{h}_k(\zeta)\to 1 $ as $k\to \infty$ uniformly for $|\zeta|\leq n_k^2 r_k$ with $r_k\to 0$.
\end{lemma}
\begin{proof}
    We have $G_k(0) =n_k\theta_{n_k}(a)=\Re\tau_{n_k}$ and
    \begin{align*}
        G_k(z) -G_k(0)&= \int_0^z G_k'(\zeta)d\zeta = \int_0^z \frac{-in_k}{\pi(n_kq)^2}m_{E_{n_k}}\left(a+\frac{\zeta}{(n_kq)^2}\right) d\zeta.
    \end{align*}
    Now, define
    \begin{align*}
        \Tilde{h}_k(\zeta) = -\frac{1}{\pi q}m_{E_{n_k}}\left(a+\frac{\zeta}{(n_kq)^2}\right) \sqrt{\frac{\zeta}{(n_kq)^2 }}\frac{\sqrt{a+\frac{\zeta}{(n_kq)^2 }-u_n}\sqrt{a+\frac{\zeta}{(n_kq)^2 }-v_n}}{a+\frac{\zeta}{(n_kq)^2 }-\lambda_n}.
    \end{align*}
    We get
    \begin{align*}
        G_k(z) -G_k(0)&= \int_0^z \frac{in_kq}{(n_kq)^2} \Tilde{h}_k(\zeta)\frac{a+\frac{\zeta}{(n_kq)^2 }-\lambda_n}{\sqrt{\frac{\zeta}{(n_kq)^2 }}\sqrt{a+\frac{\zeta}{(n_kq)^2 }-u_n}\sqrt{a+\frac{\zeta}{(n_kq)^2 }-v_n}}d\zeta\\
        &= i\int_0^z \Tilde{h}_k(\zeta) \frac{\zeta- \Tilde{\lambda}_k}{\sqrt{\zeta(\zeta-\Tilde{u}_k)(\zeta-\Tilde{v}_k)}}d\zeta.
    \end{align*}
    If $|\zeta|\leq n_k^2 r_k$, then
    \begin{align*}
        \left|a+\frac{\zeta}{(n_kq)^2} - a\right| \leq \frac{r_k}{q^2}.
    \end{align*}
    and by an application of Lemma \ref{lem:6} similar as in the proof of Proposition \ref{prop:1}
    \begin{align*}
        \Tilde{h}_k(\zeta) \xrightarrow{k\to\infty} 1,
    \end{align*}
    uniformly in $\zeta$.
\end{proof}
Next, we will need relations between the values of $\Tilde{\lambda}_k,\Tilde{u}_k$ and $ \Tilde{v}_k$ which are approximate versions of \eqref{eq:ell_int_1} and \eqref{eq:ell_int_2}.
\begin{lemma}\label{lem:12}
    \begin{multline}\label{eq_relation1}
         (\Tilde{v}_k- \Tilde{\lambda}_k)\int_0^{\pi/2} \Tilde{h}_k\left(\Tilde{u}_k \sin^2(\theta)\right) \frac{1}{\sqrt{1 - \frac{\Tilde{u}_k}{\Tilde{v}_k}\sin^2(\theta)}}d\theta\\-\Tilde{v}_k\int_0^{\pi/2} \Tilde{h}_k\left(\Tilde{u}_k \sin^2(\theta)\right) \sqrt{1-\frac{\Tilde{u}_k}{\Tilde{v}_k}\sin^2(\theta)}d\theta=0.
    \end{multline}
    \begin{multline}\label{eq_relation2}
        \sqrt{\Tilde{u}_k}\int_0^{\pi/2}\Tilde{h}_k(\Tilde{u}_k+(\Tilde{v}_k-\Tilde{u}_k) \sin^2(\theta)) \sqrt{(1+(\frac{\Tilde{v}_k}{\Tilde{u}_k}-1)\sin^2(\theta))}d\theta \\
        - \frac{\Tilde{\lambda}_k}{\sqrt{\Tilde{u}_k}}\int_0^{\pi/2} \frac{\Tilde{h}_k(\Tilde{u}_k+(\Tilde{v}_k-\Tilde{u}_k) \sin^2(\theta))}{\sqrt{(1+(\frac{\Tilde{v}_k}{\Tilde{u}_k}-1)\sin^2(\theta))}}d\theta = \frac{1}{2}.
    \end{multline}
\end{lemma}
\begin{proof}
    The equation \eqref{eq_relation1} follows from $G_k(\Tilde{u}_k)=G_k(0)$ combined with \eqref{eq_christoffel_schwarz} and the substitution $\zeta=\Tilde{u}_k\sin^2 (\theta)$. The interval $[u_{n_k},v_{n_k}]$ has equilibrium measure $1/{n_k}$ which translates to \[
    G_k(\Tilde{v}_k)-G_k(\Tilde{u}_k)= n_k\Big(\theta_{n_k}(v_{n_k})-\theta_{n_k}(u_{n_k})\Big)=n_k\mu_{\sE_{n_k}}([u_{n_k},v_{n_k}])=1.\] Then, using \eqref{eq_christoffel_schwarz} and the substitution $s=\Tilde{u}_k+(\Tilde{v}_k-\Tilde{u}_k) \sin^2(\theta)$ yields \eqref{eq_relation2}.  
\end{proof}
The following corollary can be proved in the same way as Lemma \ref{lem:18}.
\begin{corollary}\label{cor:2}
    Suppose that
    \begin{align*}
    \lambda_0 = \lim_{k\to \infty} \Tilde{\lambda}_k\in [0,\infty],\\
    u_0 = \lim_{k\to \infty} \Tilde{u}_k\in [0,\infty],\\
    v_0 = \lim_{k\to \infty} \Tilde{v}_k\in [0,\infty].
    \end{align*}
    Then, if $v_0=\infty$, then $\lambda_0=u_0=v_0=\infty$ and $\lim_{k\to\infty}\frac{\Tilde{\lambda}_k}{\Tilde{u}_k}=\lim_{k\to\infty}\frac{\Tilde{u}_k}{\Tilde{v}_k}=1$. If $0<\lambda_0<\infty$, then $\lambda_0<u_0<v_0$ and they satisfy \eqref{eq:ell_int_1} and \eqref{eq:ell_int_2}.
\end{corollary}
\subsection{Small shrinking gap}
In the following, we will further distinguish different cases depending on the value of $\lambda_0=\lim_{k\to\infty}(n_kq)^2(\lambda_{n_k}-a)$. Let us start with $\lambda_0=0$.
\begin{proposition}\label{prop:2}
    Let $n_k$ be a subsequence such that $\tau_{n_k}\to \tau\in \Upsilon$, $u_{n_k}>a$ and $u_{n_k}\to a$ as $k\to\infty$. If $(n_kq)^2(\lambda_{n_k}-a)\xrightarrow{k\to \infty} 0$, then $\tau\in\bbZ$ and
    \begin{align*}
        T_{n_k}\left(a+\frac{z}{(n_kq)^2}\right)\xrightarrow{k\to\infty}\cos\left(\pi\Theta(z,\tau)\right),
    \end{align*}
    locally uniformly for $z\in \bbC$.
\end{proposition}
\begin{proof}
    Let $G_k, \Tilde{\lambda}_k, \Tilde{u}_k,\Tilde{v}_k$ be defined by \eqref{eq:def_g_k}--\eqref{eq:tildeV}. From \eqref{eq:bound_imag_nu}, we know
    \[
    0<\Im \tau_{n_k} \xrightarrow{k\to \infty}0.
    \]
    Thus, $\tau \in \bbZ$. Since \[T_{n_k}\left(a+\frac{z}{(n_kq)^2}\right)=\cos (\pi G_k(z)),\] it remains to show
    \[
    G_k(z) \xrightarrow{k\to \infty} \tau + 1 +2i\sqrt{z-\frac{1}{4}}
    \]
    for $z\in \bbC_+$. 
    
    Let us consider $F_n = \sE_n \cup [a,u_n]$ and the comb map $\Tilde{\theta}_n$ with $\Im \Tilde{\theta}_n = \frac{1}{\pi}g_{F_n}$ and $\Tilde{\theta}_n(v_n)=0$. Since $g_{\sE_n}-g_{F_n} $ is harmonic in $ \bC\setminus F_n$ the maximum principle shows 
    \begin{align*}
        0&\leq \Im n_k\theta_{n_k}(z) -\Im n_k\Tilde{\theta}_{n_k}(z)\\
        &= \frac{n_k}{\pi}\Big(g_{\sE_{n_k}}(z)-g_{F_{n_k}}(z)\Big) \\
        &\leq  \frac{1}{\pi}n_kg_{\sE_{n_k}}(\lambda_{n_k})=\Im \tau_{n_k}\xrightarrow{k\to\infty} 0.
    \end{align*}
    This means, that 
    \begin{align}\label{eq:diff_thetas}
        n_k\theta_{n_k}\left(v_{n_k} + \frac{z}{(n_kq)^2}\right) -n_k\Tilde{\theta}_{n_k}\left(v_{n_k} + \frac{z}{(n_kq)^2}\right)
    \end{align}
    is a Herglotz function which (by possibly passing to another subsequence) converges to a (real) constant as $k\to\infty$. This constant cannot be $\infty$ since both terms are locally bounded for $z\in \bbC_+$ with $\Re z<-\Tilde{v}_k$. Indeed, the formula \eqref{eq_christoffel_schwarz} shows $\limsup_{k\to \infty}|G_k(z)|\leq2\pi\sqrt{|z|}$ for $z\in \bbC_+$ with $\Re z<0$ and we have
    \[
    n_k\theta_{n_k}\left(v_{n_k} + \frac{z}{(n_kq)^2}\right)=G_k(z+\Tilde{v}_k),
    \]
    with $\limsup_{k\to \infty}\Tilde{v}_k<\infty$ by Corollary \ref{cor:2}. Furthermore, we have seen above that $n_k(g_{\sE_{n_k}}(z)-g_{F_{n_k}}(z))\xrightarrow{k\to \infty}0$ uniformly for $z\in \bbC$. With Lemma \ref{lem:14}, this implies $g_{F_n}(z)\to g_{\sE}(z)$ locally uniformly for $z\in \bbC\setminus [\min\sE,\max\sE]$. Then, Lemma \ref{lem:13} shows that $K_n=F_n$ satisfies the assumptions of Proposition \ref{prop:1} (ii) and together with $\Tilde{\theta}_{n_k}\left(v_{n_k}\right)=0$ we get
    \begin{equation}\label{eq:theta_tilde_lim}
    n_k\Tilde{\theta}_{n_k}\left(v_{n_k} + \frac{z}{(n_kq)^2}\right)\xrightarrow{k\to \infty}2 i\sqrt{z}.
    \end{equation}
    Therefore, the limit in \eqref{eq:diff_thetas} is constant and finite. In fact, the constant is $\tau +1$ since $n_k\theta_{n_k}(v_{n_k})=\Re\tau_{n_k} +1= \tau +1$ for $k$ large enough and $n_k\Tilde{\theta}_{n_k}\left(v_{n_k}\right)=0$. Moreover, from \eqref{eq:theta_tilde_lim} we obtain
    \[
    G_k(z+\Tilde{v}_k)= n_k\theta_{n_k}\left(v_{n_k} + \frac{z}{(n_kq)^2}\right) \xrightarrow{k\to \infty} \tau + 1 +2 i\sqrt{z}.
    \]
    Now, $G_k(0)= \Re \tau_{n_k}=\tau$ for $k$ large enough which implies $\Tilde{v}_k \xrightarrow{k\to \infty}\frac{1}{4}$ and
    \[
    G_k(z)\xrightarrow{k\to \infty}\tau +1  +2 i\sqrt{z-\frac{1}{4}}.
    \]
    This concludes the proof.
\end{proof}
\subsection{Medium shrinking gap}
In this section we discuss the case $\lambda_0\in(0,\infty).$ This is precisely the case which is responsible for the appearance of the elliptic integrals introduced in Section 3. 
\begin{proposition}\label{prop:3}
    Let $n_k$ be a subsequence such that $\tau_{n_k}\to \tau\in \Upsilon$, $u_{n_k}>a$ and $u_{n_k}\to a$ as $k\to\infty$. If $(n_kq)^2(\lambda_{n_k}-a)\xrightarrow{k\to \infty} \lambda_0\in (0,\infty)$, then $\tau=m+iy$ for $m\in\bbZ$, $y>0$ and 
    \begin{align*}
        T_{n_k}\left(a+\frac{z}{(n_kq)^2}\right)\xrightarrow{k\to\infty}\cos\left(\pi\Theta(z,\tau)\right),
    \end{align*}
    locally uniformly for $z\in \bbC$.
\end{proposition}
\begin{proof}
    Let $G_k$ be the rescaled comb map \eqref{eq:def_g_k} and $\lambda_0,u_0,v_0$ the limits of $\Tilde{\lambda}_k,\Tilde{u}_k,\Tilde{v}_k$ respectively. We know from Corollary \ref{cor:2} that $0<\lambda_0<u_0<v_0<\infty$ and that they satisfy \eqref{eq:ell_int_1}, \eqref{eq:ell_int_2}. Then, bounded convergence and \eqref{eq_christoffel_schwarz} shows that
    \[
    G_k(z) \xrightarrow{k\to \infty} \Re \tau + i\int_0^z \frac{\zeta-\lambda_0}{\sqrt{\zeta(\zeta-u_0)(\zeta-v_0)}}d\zeta
    \]
    locally uniformly for $z\in \bbC_+$. Let us denote the limit as $\Theta_0$. By Lemma \ref{lem:17}, we know that $\Theta_0=\Theta(\cdot,\Re\tau+iy)$ for some $y>0$. 
    Moreover, \eqref{eq_christoffel_schwarz} and bounded convergence show that for a fixed $z\in \bbC_+$
    \[
    G_k(z)- G_k(x) \xrightarrow{k\to \infty}i\int_x^z \frac{\zeta-\lambda_0}{\sqrt{\zeta(\zeta-u_0)(\zeta-v_0)}}d\zeta= \Theta_0(z)-\Theta_0(x)
    \]
    pointwise for $x\in \bbR$ and uniformly for $x$ in compact subsets of $\bbR\setminus \{0,u_0,v_0\}$. This implies $G_k\to \Theta_0$ uniformly on compact subsets of $\bbC_+\cup \bbR$. In particular, since $0<\lambda_0<u_0$, 
    \[
    \Re \tau +iy=\Theta_0(\lambda_0) =\lim_{k\to \infty}G_k(\Tilde{\lambda}_k)=\lim_{k\to\infty}\tau_{n_k}=\tau.
    \]
    Therefore, $\Theta_0=  \Theta(\cdot,\tau)$ and
    \[
    T_{n_k}\left(a+\frac{z}{(n_kq)^2}\right)=\cos(\pi G_k(z))\xrightarrow{k\to\infty}\cos\left(\pi\Theta(z,\tau)\right)
    \]
    for $z\in \bbC_+$. By Corollary \ref{cor:3} the convergence extends to locally uniform convergence in $\bbC$.
\end{proof}
\subsection{Large shrinking gap}
\begin{proposition}\label{prop:4}
    Let $n_k$ be a subsequence such that $\tau_{n_k}\to \tau\in \Upsilon$, $u_{n_k}>a$ and $u_{n_k}\to a$ as $k\to\infty$. If $(n_kq)^2(\lambda_{n_k}-a)\xrightarrow{k\to \infty} \infty$, then $\tau=m+i\infty$ for $m\in\bbZ$ and 
    \begin{align*}
        T_{n_k}\left(a+\frac{z}{(n_kq)^2}\right)\xrightarrow{k\to\infty}\cos\left(\pi\Theta(z,\tau)\right),
    \end{align*}
    locally uniformly for $z\in \bbC$.
\end{proposition}
\begin{proof}
    Let $G_k$ be defined by \eqref{eq:def_g_k}. Since $\Tilde{\lambda}_k\to \infty$, Corollary \ref{cor:2} shows that \[\lim_{k\to\infty}\frac{\Tilde{\lambda}_k}{\Tilde{u}_k}=\lim_{k\to\infty}\frac{\Tilde{u}_k}{\Tilde{v}_k}=1.\] Then, \eqref{eq_christoffel_schwarz} shows that
    \[
    G_k(z) \xrightarrow{k\to \infty} \Re \tau + 2i\sqrt{z}
    \]
    uniformly for $z$ in compact subsets of $\bbC_+\cup \bbR$. Let $M>0$. Then  for $k$ large enough, we have $2M^2 <\Tilde{\lambda}_k $ and $M<\Im G_k (2M^2)<\Im \tau_{n_k}$. Therefore, $\tau = \Re \tau + i\infty$ and 
    \[
    \Theta(z,\tau)= \Re \tau + 2i\sqrt{z}.
    \]
    This proves the assertion.
\end{proof}
\begin{prf}{Theorem \ref{mainTheorem}}
    Let us suppose 
    \begin{equation}\label{eq:diff_resc}
    T_n\left(a+\frac{z}{(nq)^2}\right)-\cos(\pi\Theta(z,\tau_n))
    \end{equation}
    does not converge to $0$. Then, we can find a subsequence $n_k$ such that any further subsequence also does not converge to $0$ locally uniformly in $\bbC$. Since we might assume $\Re \tau_{n_k}\in [0,2]$, we can suppose that $\tau_{n_k}\to \tau\in \Upsilon$. Then, Proposition \ref{prop:7} shows
    \[
    \cos(\pi\Theta(z,\tau_{n_k}))\xrightarrow{k\to \infty}\cos(\pi\Theta(z,\tau))
    \]
    locally uniformly for $z\in \bbC$. Furthermore, taking a further subsequence, we can assume that $n_k$ satisfies the assumptions of Proposition \ref{prop:2}, \ref{prop:3}, \ref{prop:4}, \ref{prop:5} or \ref{prop:6}. Then,
    \[
    T_{n_k}\left(a+\frac{z}{(n_kq)^2}\right)\xrightarrow{k\to \infty} \cos (\pi\Theta(z,\tau))
    \]
    locally uniformly for $z\in \bbC$. However, this contradicts 
    \[
    T_{n_k}\left(a+\frac{z}{(n_kq)^2}\right)-\cos(\pi\Theta(z,\tau_{n_k}))\not\to 0.
    \]
    Thus, \eqref{eq:diff_resc} must converge to $0$.
\end{prf}
\begin{prf}{Corollary \ref{cor:1}}
    Let $x_j^{(n)}$ be the zeros of $T_n$ relative to $a$. Furthermore, let $x_j^\tau$ be the zeros of $\cos(\pi \Theta(\cdot,\tau))$ for $\tau\in \Upsilon$ as in Lemma \ref{lem:15}. Let us suppose that there exists a subsequence $n_k$ such that $n_k^2(x_{l}^{(n_k)}-x_j^{(n_k)})$ converges to either $0$ or $\infty$ for $j<l$. By taking a further subsequence we can suppose $\tau_{n_k}\to \tau\in \Upsilon$. Then, Theorem \ref{mainTheorem} and Hurwitz's theorem show 
    \[
    (n_kq)^2 (x_j^{(n_k)}-a) \xrightarrow{k\to \infty} x_j^\tau,\\
    (n_kq)^2 (x_l^{(n_k)}-a) \xrightarrow{k\to \infty} x_l^\tau,
    \]
    except when $\tau_{n_k}\to (m+\frac{1}{2})^-$ where the limits are $x_{j-1}^\tau$ and $x_{l-1}^\tau$. In any case, the difference between those two limits is bounded above and below by positive constants by Lemma \ref{lem:15}. This contradicts $\lim_{k\to \infty}n_k^2(x_{l}^{(n_k)}-x_j^{(n_k)})\in \{0,\infty\}$. Since $\cos(\pi\Theta(\cdot,\tau))$ has at most one zero with positive real part, one can show that $T_n$ can only have at most one zero in $[a,a+Rn^{-2}]$ for every $R>0$ and $n$ large enough. Note that the above argument also works for $l=1$ when we restrict to subsequences where $x_1^{(n_k)}$ exists. 
\end{prf}
\section{Generic Akhiezer polynomials accumulate in all of $\Upsilon$}
In this section, we will focus on the case of two intervals, already studied by Akhiezer \cite{Akhiezer1,Akhiezer2,Akhiezer3}. Let us recall some notation. For $-1<a<b<1$ and $\sE=[-1,a]\cup [b,1]$, the $n$-extension $\sE_n=(T_n^\sE)^{-1}([-1,1])$ can only have one extra interval $[u_n,v_n]$ in the gap $[a,b]$ with the convention $u_n=v_n=b$ if there is no extra interval. Then, $\lambda_n$ is the critical point of $T_n^\sE$ in $[a,u_n]$ if $u_n>a$ and $\lambda_n =a$ otherwise. We consider the parameter $\tau_n = \frac{1}{\pi}\cos^{-1}(T_n(\lambda_n))\in \Upsilon$. Since we will vary $a,b$, it will be convenient to consider $u_n(a,b), v_n(a,b),\lambda_n(a,b), \tau_n(a,b)$ as functions of $a,b$. To prove Theorem \ref{thm:all_limits} our goal is to find a set $F$ such that for all $(a;b)\in F$ the sequence $\tau_n(a,b)$ accumulates in all of $\Upsilon$ (modulo adding real part $2k$) as $n\to\infty$. We start by showing continuity of $\tau_n(a,b)$. To exclude some edge cases, let $\Upsilon^* = \Upsilon \setminus \{k+iy: k\in \bbZ, y\in \{0,\infty\}\}$. 
\begin{lemma}\label{lem:nu_cont}
    Let $n\in \bbN$. Then, \[B_n :=\tau_n^{-1}(\Upsilon^*)\cap \{(a;b): v_n(a,b)<b\}\] is open in $(-1,1)^2$ and $\tau_n$ restricted to $B_n$ is continuous. 
\end{lemma}
\begin{proof}
    Let $-1<a<b<1$ with $(a;b)\in B_n$ and $a_k,b_k$ two sequences with $-1<a_k<b_k<1 $ and $a_k\xrightarrow{k\to \infty}a,b_k\xrightarrow{k\to \infty}b$. We have to show that $(a_k;b_k)\in B_n$ for big enough $k$ and $\tau_n(a_k,b_k)\xrightarrow{k\to\infty} \tau_n(a,b)$. Define $\sE=[-1,a]\cup [b,1]$ and $\sE_k=[-1,a_k]\cup[b_k,1]$. \cite[Satz III]{Akhiezer1} shows that the Chebyshev polynomial $T_n^{\sE_k}$ converges to $T_n^\sE$ locally uniformly in $\bbC$ as $k\to\infty$. Since $v_n(a,b)<b$ we know there are two possible cases, either $a=u_n(a,b)$ or $a<u_n(a,b)$. We will consider each case individually. 
    
    First assume $a=u_n(a,b)$. Since $T_n^{\sE_k}$ converges locally uniformly to $T_n^\sE$, we know $T_n^{\sE_k}(a_k)\xrightarrow{k\to\infty}T_n^\sE(a)$. Since $\tau_n(a,b)\in\Upsilon^*$ and thus $\tau_n(a,b)\notin \bbZ$, we have $|T_n^\sE(a)|<1$ and thus for big enough $k$ we have $|T_n^{\sE_k}(a_k)|<1$. This implies $a_k=u_n(a_k,b_k)$ and $\lambda_n(a_k,b_k)=a_k$ and we get $T_n^{\sE_k}(\lambda_n(a_k,b_k))\xrightarrow{k\to\infty} T_n^\sE(a)$. Thus, $\tau_n(a_k,b_k)\xrightarrow{k\to \infty}\tau_n(a,b)$. Since $|T_n^\sE(x)|>1$ for $v_n(a,b)<x<b$, we also have $|T_n^{\sE_k}(x)|>1$ for some $x\in (a_k,b_k)$ and $k$ big enough. Therefore, $v_n(a_k,b_k)<b_k$ and $|T_n^{\sE_k}(a_k)|<1$ for big enough $k$, which implies $(a_k;b_k)\in B_n$.

    Now, let $a<u_n(a,b)$. Since $v_n(a,b)<b$, we also have $u_n(a,b)<v_n(a,b)$. Then, $T_n^\sE$ has a simple zero $ \xi$ in $(a,b)$. For big enough $k$, $T_n^{\sE_k}$ must have a simple zero $\xi_k\in (a_k,b_k)$ with $\xi_k\to \xi$.  As before, we also have $|T_n^\sE(x)|>1$ for $a<x<u_n(a,b) $ and $v_n(a,b)<x<b$. This is easily seen to give $a_k<u_n(a_k,b_k)<v_n(a_k,b_k)<b$ for big enough $k$, which means $(a_k;b_k)\in B_n$. We also know that $\lambda_n(a,b)$ is the unique simple zero of $(T_n^\sE)'$ in the interval $(a,\xi)$. Since $\lambda_n(a_k,b_k)$ is the zero of $(T_n^{\sE_k})'$ in $(a_k,\xi_k)$, this gives $\lambda_n(a_k,b_k)\xrightarrow{k\to\infty} \lambda_n(a,b)$. We get $\tau_n(a_k,b_k)\xrightarrow{k\to\infty}\tau_n(a,b)$, proving the lemma.
\end{proof}
\begin{figure}
    \centering
    	\begin{tikzpicture}[scale=0.7]
		\begin{scope}[xshift=-8cm]
			\draw[-, thick] (-6,0) -- (-6,5);
            \draw[-, thick] (0,0) -- (-6,0);
		      \draw[-, thick] (-2,0) -- (-2,3);
            \draw[-, thick] (-2.5,0) -- (-2.5,1);
		      \draw[-, thick] (0,0) -- (0,5);

            \fill (-2.5,1) circle (3pt);
            \fill (-2,3) circle (3pt);
            \node at (-0,-0.5) {$0$};
            \node at (-5.5,-0.5) {$-1$};
            \node at (-2.1,3.5) {$\xi_n+\frac{1}{n}+ih$};
            \node at (-3.4,1.5) {$\xi_n+i\frac{\Im \tau}{n}$};
		
		\end{scope}
		
		\draw[-, thick] (-6,0) -- (-6,5);
        \draw[-, thick] (0,0) -- (-6,0);
		\draw[-, thick] (-2,0) -- (-2,3);
		\draw[-, thick] (0,0) -- (0,5);

            \fill (-2,3) circle (3pt);
            \node at (-0,-0.5) {$0$};
            \node at (-5.5,-0.5) {$-1$};
            \node at (-1.5,3.5) {$\xi+ih$};
		
	\end{tikzpicture}
    \caption{Approximation of the comb $\Lambda$ by $\Lambda_n$ in Lemma \ref{lem:9}.}
    \label{fig:comb}
\end{figure}
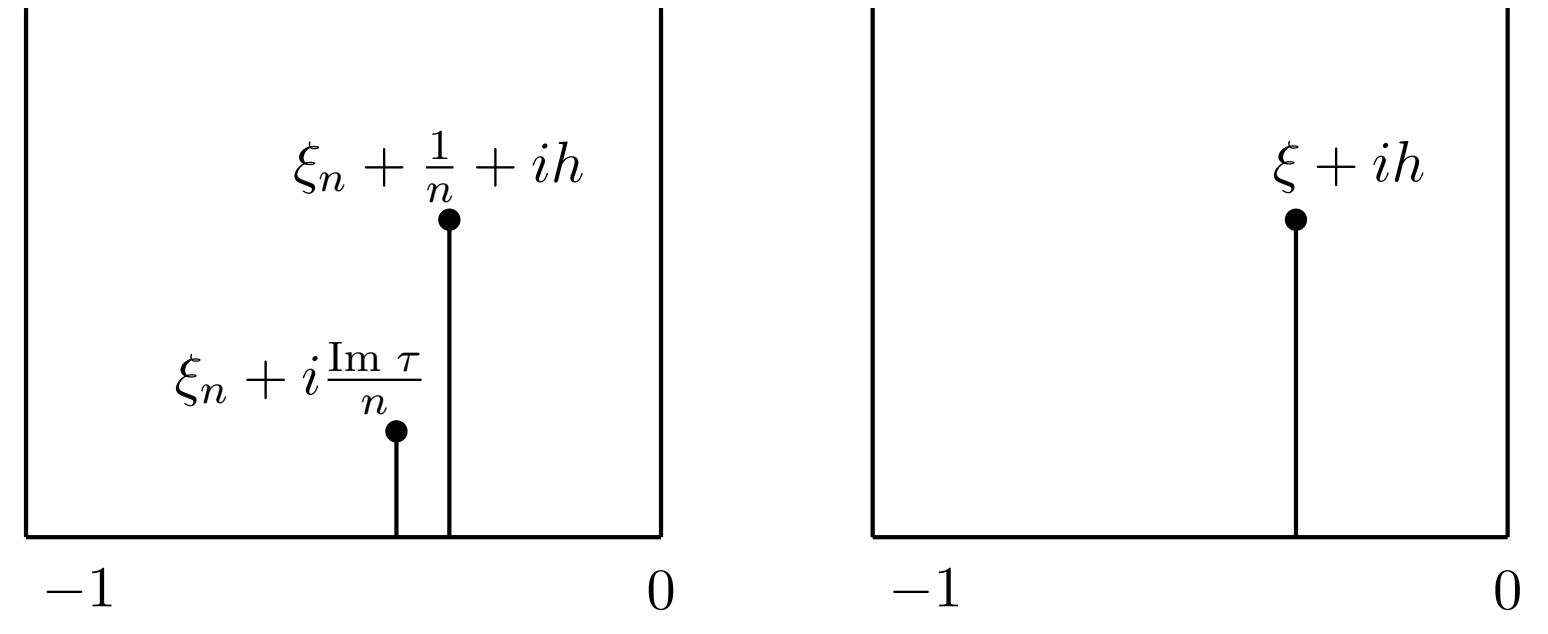
\begin{lemma}\label{lem:9}
    Let $-1<a<b<1$ and $\tau\in \Upsilon^*$ with $0\leq \Re \tau<2 $. Then there exists a sequence $(a_n,b_n)$ such that $(a_n;b_n)\in B_n$, $a_n\to a$, $b_n\to b$ and $\tau_n(a_n,b_n)=\tau+2m_n$ for some $m_n\in \bbZ$.
\end{lemma}
\begin{proof}
    Let $\sE = [-1,a]\cup [b,1]$. Then, we consider the comb map $\theta:\bbC_+\to \Lambda$ such that $\Im \theta=\frac{1}{\pi}g_\sE$ and
    \begin{align*}
        \Lambda=\{z\in \bbC_+: -1 < \Re z < 0\} \setminus\{\xi + i y: 0<y\leq h\},
    \end{align*}
    for some $-1<\xi<0$, $h >0$. Note that $\theta(-1)=-1, \theta(1)=0$. We will construct a sequence of combs $\Lambda_n$ that approximate the comb $\Lambda$, see Figure \ref{fig:comb}. Let $m_n = \left\lfloor\frac{n}{2}\xi \right\rfloor$, $\xi_n =\frac{2m_n+ \lfloor\Re \tau\rfloor}{n}$ and
    \begin{align*}
        \mathcal{N}_n &= \left\{ \xi_n + iy : 0<y\leq \frac{\Im \tau}{n}\right\}\cup\left\{ \xi_n+\frac{1}{n} + iy : 0<y\leq h\right\}.
    \end{align*}
    We define $\Lambda_n= \{z\in \bbC_+: -1 < \Re z < 0\} \setminus\mathcal{N}_n$ and let $\theta_n$ be the conformal map from $\bbC_+$ to $\Lambda_n$ with $\theta_n(-1)=-1, \theta_n(1)=0$ and $\theta_n(\infty)=\infty$. 
    
    Then, we set \[a_n = \theta_n^{-1}\left(\frac{2m_n +\Re \tau}{n}-0\right), \quad \text{and}\quad b_n =  \theta_n^{-1}\left(\xi_n+\frac{1}{n}+0\right).\] The $n$-th Chebyshev polynomial of $[-1,a_n]\cup [b_n,1]$ is given by $\cos (n\pi\theta_n)$. We consider $u_n = u_n(a_n,b_n), v_n=v_n(a_n,b_n)$ and $\lambda_n =\lambda_n(a_n,b_n)$. Note that $\tau_n(a_n,b_n)=n \left(\frac{2m_n +\Re \tau}{n} + i\frac{\Im \tau}{n} \right)= \tau + 2m_n$. Since $\tau \in\Upsilon^*$, we also get $\tau_n(a_n,b_n)\in \Upsilon^*$. From the comb picture, it is also clear that $v_n(a_n,b_n)<b_n$, i.e. $(a_n;b_n)\in B_n$. It remains to show that $a_n\to a$ and $b_n\to b$. 

    We know, $\sE_n=[-1,a_n]\cup [u_n,v_n]\cup [b_n,1]$ is the $n$-extension of $[-1,a_n]\cup [b_n,1]$. Similarly as in Lemma \ref{lem:6} one can show that
    \begin{align*}
        m_{\sE_n}(z) = \int_\bbR \frac{1}{t-z}d\mu_{\sE_n}(t) = - \frac{(z-\lambda_n)(z-\Tilde{\lambda}_n)}{\sqrt{(z^2-1)(z-a_n)(z-b_n)(z-u_n)(z-v_n)}},
    \end{align*}
    where $\Tilde{\lambda}_n$ is the critical point of $g_{\sE_n}$ in $(v_n,b_n)$. We also know $\theta_n(z)=-\frac{i}{\pi}\int_0^z m_{\sE_n}(t)dt$. 
    
    Let us now consider a subsequence $n_k$ such that $a_{n_k},b_{n_k},\lambda_{n_k},\Tilde{\lambda}_{n_k},u_{n_k},v_{n_k}$ converge to $a_0,b_0,\lambda,\Tilde{\lambda},u,v$ respectively. Then,
    \begin{align*}
        m_{\sE_{n_k}}(\zeta) \xrightarrow{k\to\infty} - \frac{(\zeta-\lambda)(\zeta-\Tilde{\lambda})}{\sqrt{(\zeta^2-1)(\zeta-a_0)(\zeta-b_0)(\zeta-u)(\zeta-v)}}
    \end{align*}
    uniformly for $\zeta$ in compact subsets of $\bbC_+\cup \bbR \setminus\{-1,1,a_0,b_0,u,v\}$. Let us denote the RHS with $f(\zeta)$. Then, we get
    \begin{align}\label{eq:estimate_f}
        \theta_{n_k}(z)-\theta_{n_k}(w)\xrightarrow{k\to\infty} -\frac{i}{\pi}\int_w^zf(\zeta)d\zeta,
    \end{align}
    uniformly for $z,w\in K$, for every convex compact subset $K\subset\bbC\cup \bbR\setminus \{-1,1,a_0,b_0,u,v\}$. 
    It is clear that if $\theta_{n_k}(x_0)$ converges for any point $x_0\in \bbC_+\cup \bbR \setminus\{-1,1,a_0,b_0,u,v\}$, then $\theta_{n_k}$ converges locally uniformly in $\bbC_+\cup \bbR \setminus\{-1,1,a_0,b_0,u,v\}$. Since $|\theta_n(x_0)|\leq 1 +h $ for $x_0\in (-1,1)$, we can always choose a further subsequence (denoted the same) such that $\theta_{n_k}\to \theta_0$ locally uniformly in $\bbC_+\cup \bbR \setminus\{-1,1,a_0,b_0,u,v\}$. Since $\theta_{n_k}$ is conformal, $\theta_0$ is either conformal or constant. The latter is not possible because $\theta_0'(z)=f(z)\ne 0$. By \eqref{eq:estimate_f}, we get $\int_I|f(\zeta)|d\zeta\leq 2+2h$ for every interval $I\subseteq (-1,1)\setminus \{a_0,b_0,u,v\}$. This implies that $f$ is integrable over $[-1,1]$ and can only have square-root singularities. Then, $\theta_0$ extends to $\bbR$ and is a Christoffel-Schwarz mapping. 
    Moreover, we know $-1\leq a_0\leq \lambda\leq u\leq v\leq \Tilde{\lambda}\leq b_0\leq 1$ and using $\theta_{n_k}\to \theta_0$ on $\bbR\setminus \{-1,1,a_0,b_0,u,v\}$ we get
    \begin{align*}
        \theta_0((-\infty,-1))&\subseteq (-1+i(\bbR_{\geq 0})\cup \{\infty\},\\
        \theta_0((-1,a_0))&\subseteq [-1,\xi],\\
        \theta_0((a_0,v))&\subseteq \{\xi\},\\
        \theta_0((v, b_0))&\subseteq \{\xi+iy : 0\leq y\leq h\},\\
        \theta_0((b_0,1))&\subseteq [\xi,0],\\
        \theta_0((1,\infty))&\subseteq i\bbR_{\geq 0}\cup \{\infty\}.
    \end{align*}
    We can draw several conclusions from this. Firstly, $\theta_0(-1)=-1, \theta_0(1)=0,\theta_0(\infty)=\infty$ and $\theta_0(a_0)=\theta_0(b_0)=\xi$. Furthermore, we also know $a_0=\lambda=u=v$ and that $\theta_0([v,b_0])=\{\xi+iy : 0\leq y\leq \Tilde{h}\}$ where $\theta_0(\Tilde{\lambda})=m+i\Tilde{h}$ and $\Tilde{h}\leq h$. This means that
    \[
    \theta_0(\bbC_+)=\{z\in \bbC_+:-1<\Re z<0\}\setminus \{\xi+iy:0<y\leq \Tilde{h}\}.
    \]
    To show $\Tilde{h}=h$, we follow an argument also given in \cite[Proof of Theorem 1.8]{Pom92}. If $\Tilde{h}<h$, then we can find $z=\xi + iy$ with $\Tilde{h}<y<h$ and $\delta>0$ such that
    \[
    A := \{w\in \bbC: |w-z|\leq 2\delta\}\subseteq \theta_0(\bbC_+).
    \]
    Since $\theta_0$ is a conformal map, we know that $B=\theta_0^{-1}(A)\subseteq\bbC_+$ is compact. Then, we can find $k_0$ such that $|\theta_{n_k}(w)-\theta_0(w)|<\delta$ for all $k\geq k_0$ and $w\in B$. For every $w_0\in \bbC$ with $ |w_0-z|<\delta$, we know $|\theta_0(w)-w_0|\geq |\theta_0(w)-z|-\delta=\delta$ for $w\in \partial B$. Since $\theta_0-w_0$ has a zero in $B$, Rouch\'e's theorem implies that $\theta_{n_k}-w_0$ also has a zero in $B$ for $k\geq k_0$, i.e.
    \[
        \{w\in \bbC: |w-z|<\delta\}\subseteq \theta_{n_k}(\bbC_+).
    \]
    However, this cannot hold true since
    \[
    w_{n_k}=\frac{2m_{n_k}+\lfloor \Re \tau\rfloor}{n_k}+iy \in \mathcal{N}_{n_k} \subseteq \bbC_+\setminus \theta_{n_k}(\bbC_+)
    \]
    but $w_{n_k}\to z$. Thus, we have $\Tilde{h}=h$ and $\theta_0(\bbC_+)=\Lambda$. Then, $\theta_0 = \theta$ and $a_0=\theta_0^{-1}(\xi-)=\theta^{-1}(\xi-)=a$, as well as $b_0=b$ which concludes the proof.
\end{proof}
Now, we can prove Theorem \ref{thm:all_limits}.
\begin{prf}{Theorem \ref{thm:all_limits}}
    For $\tau\in \Upsilon^*$ let 
    \begin{align*}
        F_n(\tau) : = B_n \cap \tau_n^{-1} \left(\bigcup_{m\in \bbZ}\{z\in \bbC : |z-(\tau+2m)|< \frac{1}{n}\}\right).
    \end{align*}
    By Lemma \ref{lem:nu_cont}, we know that $F_n(\tau)$ is open. By Lemma \ref{lem:9} the set
    \begin{align*}
        \bigcup_{n\geq n_0} F_n(\tau)
    \end{align*}
is dense in $\{(a;b): -1 <a<b<1\}$ for every $n_0\in \bbN$. Indeed, for every $-1<a<b<1$, we find $(a_n;b_n)\in B_n$ with $\tau_n(a_n,b_n)=\tau+2m_n$ (i.e. $(a_n;b_n)\in F_n(\tau)$) and $a_n \to a, b_n \to b$. Let $\bbQ\Upsilon^* = \{\tau\in \Upsilon^*: \Re\tau,\Im \tau \in \bbQ\}$. Then, by Baire's category theorem
\begin{align*}
    F := \bigcap_{\tau \in\bbQ\Upsilon^*}\bigcap_{n_0\in \bbN}\bigcup_{n\geq n_0} F_n(\tau)
\end{align*}
is a dense $G_\delta$ set in $\{(a;b): -1 <a<b<1\}$. For every $(a;b)\in F$ and $ \tau \in \bbQ\Upsilon^*$, we find a subsequence $n_k$ such that $\tau_{n_k}(a,b)-2m_k\xrightarrow{k\to\infty}\tau$. Proposition \ref{prop:7} shows that  $\cos(\pi\Theta(\cdot,\tau_{n_k}(a,b)))=\cos(\pi\Theta(\cdot,\tau_{n_k}(a,b)-2m_k))\xrightarrow{k\to\infty}\cos(\pi\Theta(\cdot,\tau))$ locally uniformly in $\bbC$. Combined with Theorem \ref{mainTheorem}, this implies $T_{n_k}^\sE(a+ \frac{z}{(n_kq)^2})\xrightarrow{k\to \infty} \cos(\pi \Theta(\cdot,\tau))$ for $\sE=[-1,a]\cup [b,1]$ and locally uniformly in $\bbC$. Thus, $\cos(\pi\Theta(\cdot,\tau))$ is a limit point of the sequence $T_n^{\sE}\left(a+ \frac{\cdot}{(nq)^2}\right)$ for every $\tau \in \bbQ\Upsilon^*$. However, the set of limit points is closed and $\bbQ\Upsilon^*$ is dense in $\Upsilon$ which means that the same holds for all $\tau\in \Upsilon$.
\end{prf}

\section*{Acknowledgements} B.B. was supported by Grant 1154426N from the Research Foundation
Flanders (FWO). B.B. thanks Athanasios Kouroupis and Joakim Cronvall for a helpful lecture on Carath\'eodory kernel convergence and boundary behavior of conformal maps.  B.E. is grateful to Peter Yuditskii for a helpful discussion related to the proof of Theorem \ref{thm:all_limits}. 

\section*{AI disclosure}
Claude Sonnet 5.5 was used for text editing in the introduction and for the creation of the figures. Microsoft Copilot was used for proof reading and suggested to use elliptic integrals to represent the conditions \eqref{eq:ell_int_1} and \eqref{eq:ell_int_2} in terms of complete elliptic integrals. All other proofs and mathematical ideas are entirely due to the authors.

\providecommand{\MR}[1]{}
\providecommand{\bysame}{\leavevmode\hbox to3em{\hrulefill}\thinspace}
\providecommand{\MR}{\relax\ifhmode\unskip\space\fi MR }
\providecommand{\MRhref}[2]{%
  \href{http://www.ams.org/mathscinet-getitem?mr=#1}{#2}
}
\providecommand{\href}[2]{#2}


\begin{thebibliography}{10}

\bibitem{akhiezer.levin:1960}
N.~I. Ahiezer and B.~Ja. Levin, \emph{Generalization of {S}. {N}. {B}ernstein's inequality for derivatives of entire functions}, Issledovanija po sovremennym problemam teorii funkcii kompleksnogo peremennogo, Gosudarstv. Izdat. Fiz.-Mat. Lit., Moscow, 1960, pp.~111--165. \MR{125968}

\bibitem{Akhiezer1}
N.~Akhiezer, \emph{{\"U}ber einige {F}unktionen, welche in zwei gegebenen {I}ntervallen am wenigsten von {N}ull abweichen. {T}eil {I}}, Izv. Akad. Nauk SSSR Otdel. Mat. Estestv. Nauk. (1932), no.~9, 1163--1202.

\bibitem{Akhiezer2}
\bysame, \emph{{\"Uber einige Funktionen, welche in zwei gegebenen Intervallen am wenigsten von Null abweichen. Teil II}}, Izv. Akad. Nauk SSSR Otdel. Mat. Estestv. Nauk. (1933), no.~3, 309--344.

\bibitem{Akhiezer3}
\bysame, \emph{{\"Uber einige Funktionen, welche in zwei gegebenen Intervallen am wenigsten von Null abweichen. Teil III}}, Izv. Akad. Nauk SSSR Otdel. Mat. Estestv. Nauk. (1933), no.~4, 449--536.

\bibitem{bernstein12}
S.~N. Bernstein, \emph{D\'{e}monstration du th\'{e}or\`{e}me de {W}eierstrass fond\'{e}e sur le calcul des probabilit\'{e}s}, Comm. Kharkov Math. Soc. \textbf{13} (1912), 1--2.

\bibitem{ChriSiZi1}
J.~S. Christiansen, B.~Simon, and M.~Zinchenko, \emph{Asymptotics of {C}hebyshev polynomials, {I}: subsets of {$\Bbb R$}}, Invent. Math. \textbf{208} (2017), no.~1, 217--245. \MR{3621835}

\bibitem{eichinger.lukic.kheifets.yuditskii:2026}
B.~Eichinger, A.~Kheifets, M.~Luki\'{c}, and P.~Yuditskii, \emph{New universality classes associated to fractals}, arXiv:2607.28379 (2026), 1--78.

\bibitem{eichinger.lukic:2025}
B.~Eichinger and M.~Luki\'{c}, \emph{Stahl-{T}otik regularity for continuum {S}chr\"{o}dinger operators}, Anal. PDE \textbf{18} (2025), no.~3, 591--628. \MR{4874190}

\bibitem{eichinger.lukic.simanek:2026}
B.~Eichinger, M.~Luki\'{c}, and B.~Simanek, \emph{An approach to universality using {W}eyl {$m$}-functions}, Ann. of Math. (2) \textbf{203} (2026), no.~2, 471--510. \MR{5037837}

\bibitem{eichinger.lukic.woracek:2024}
B.~Eichinger, M.~Luki\'{c}, and H.~Woracek, \emph{Necessary and sufficient conditions for universality limits}, Invent. Math. (to appear), arXiv:2409.18045 (2024), 1--78.

\bibitem{Ahlfors}
B.~Eichinger and P.~Yuditskii, \emph{The {A}hlfors problem for polynomials}, Mat. Sb. \textbf{209} (2018), no.~3, 34--66. \MR{3769214}

\bibitem{eremenko.yuditskii:2012}
A.~Eremenko and P.~Yuditskii, \emph{Comb functions}, Recent advances in orthogonal polynomials, special functions, and their applications, Contemp. Math., vol. 578, Amer. Math. Soc., Providence, RI, 2012, pp.~99--118. \MR{2964141}

\bibitem{Faber1920}
G.~Faber, \emph{\"{U}ber {T}schebyscheffsche {P}olynome}, J. Reine Angew. Math. \textbf{150} (1920), 79--106. \MR{1580974}

\bibitem{Freud:1969}
G.~Freud, \emph{Orthogonale {P}olynome}, Lehrb\"{u}cher und Monographien aus dem Gebiete der Exakten Wissenschaften, Mathematische Reihe, Band 33, Birkh\"{a}user Verlag, Basel-Stuttgart, 1969. \MR{481888}

\bibitem{kotani:2015}
S.~Kotani, \emph{Spectral problems of ergodic {S}chr\"{o}dinger operators}, Chinese J. Appl. Probab. Statist. \textbf{31} (2015), no.~6, 596--661. \MR{3494949}

\bibitem{kuijlaars.vanlessen:2002}
A.~B.~J. Kuijlaars and M.~Vanlessen, \emph{Universality for eigenvalue correlations from the modified {J}acobi unitary ensemble}, Int. Math. Res. Not. (2002), no.~30, 1575--1600. \MR{1912278}

\bibitem{lubinsky:2009}
D.~S. Lubinsky, \emph{A new approach to universality limits involving orthogonal polynomials}, Ann. of Math. (2) \textbf{170} (2009), no.~2, 915--939. \MR{2552113}

\bibitem{lubsinksyConstr}
\bysame, \emph{{Local Limits and Spacing of Zeros of Chebyshev Polynomials on Compact Subsets of $\mathbb{R}$}}, Constr. Approx. (2026).

\bibitem{lukic:book}
M.~Luki\'{c}, \emph{A first course in spectral theory}, Graduate Studies in Mathematics, vol. 226, American Mathematical Society, Providence, RI, [2022] \copyright 2022. \MR{4497147}

\bibitem{Pom92}
Ch. Pommerenke, \emph{Boundary behaviour of conformal maps}, Grundlehren der mathematischen Wissenschaften [Fundamental Principles of Mathematical Sciences], vol. 299, Springer-Verlag, Berlin, 1992. \MR{1217706}

\bibitem{Ran95}
T.~Ransford, \emph{Potential theory in the complex plane}, London Mathematical Society Student Texts, vol.~28, Cambridge University Press, Cambridge, 1995. \MR{1334766}

\bibitem{Simon:2007}
B.~Simon, \emph{Equilibrium measures and capacities in spectral theory}, Inverse Probl. Imaging \textbf{1} (2007), no.~4, 713--772. \MR{2350223}

\bibitem{sodin.yuditskii:1992}
M.~L. Sodin and P.~Yuditskii, \emph{Functions that deviate least from zero on closed subsets of the real axis}, Algebra i Analiz \textbf{4} (1992), no.~2, 1--61. \MR{1182392}

\bibitem{szego:1924}
G.~Szeg\"o, \emph{Bemerkungen zu einer {A}rbeit von {H}errn {M}. {F}ekete: \"uber die {V}erteilung der {W}urzeln bei gewissen algebraischen {G}leichungen mit ganzzahligen {K}oeffizienten}, Math. Z. \textbf{21} (1924), no.~1, 203--208. \MR{1544698}

\bibitem{widom1969}
H.~Widom, \emph{{Extremal polynomials associated with a system of curves in the complex plane}}, Advances in Math. \textbf{3} (1969), 127--232. \MR{239059}

\end{thebibliography}
\end{document}